\documentclass[11pt]{article}
\usepackage{a4}
\usepackage{graphicx}
\usepackage{parskip}

\usepackage{amsfonts}
\usepackage{natbib}
\usepackage{algorithm}
\usepackage{algorithmic}

\usepackage{latexsym}

\newtheorem{theorem}{Theorem}[section]

\newtheorem{lemma}{Lemma}[section]

\newtheorem{remark}{Remark}[section]

\newtheorem{summary}{Summary}[section]
\newtheorem{definition}{Definition}[section] 
\newtheorem{condition}{Condition}[section]

\newcommand{\R}{\mathbb{R}}
\newcommand{\Nat}{\mathbb{N}}

\newcommand{\PP} {{  \rm I\hskip-0.22em P}}

\newcommand{\EE} {{\rm I\hskip-0.48em E}}

\usepackage{natbib}

\begin{document}

\centerline {\bf \Large On the uniform convergence of empirical norms}

\centerline {\bf \Large and inner products, }

\centerline{\bf \Large with application to causal inference}

\vskip .1in
\centerline {Sara van de Geer}

\centerline {Seminar for Statistics, ETH Z\"urich}

{\bf Abstract.} Uniform convergence of empirical norms - empirical measures of
squared functions - is a topic which has received considerable attention in the literature
on empirical processes. The
results are relevant as empirical norms
occur due to symmetrization. They also play a prominent role in statistical applications.
The contraction inequality has been a main tool but recently other approaches have shown
to lead to better results in important cases.  We present an overview including 
the linear (anisotropic) case, and give
new results for inner products of functions.  Our main application will be the estimation
of the parental structure in a directed acyclic graph. As intermediate result we establish 
convergence of the least squares estimator when the model is wrong.

\section{Introduction}\label{intro.section}

Let $X_1 , \ldots , X_n$ be independent random variables with values in ${\cal X}$ and
${\cal F} $ be a class of real-valued functions on ${\cal X}$.
For a function $f : \ {\cal X} \rightarrow \R$, we denote its empirical measure by 
$ P_n f := \sum_{i=1}^n f (X_i) / n $ and its theoretical measure by $Pf   := \sum_{i=1}^n \EE f (X_i) / n$
(assuming it exists).
Furthermore, we let $\| f \|_n^2 : = P_n f^2   $ and $ \| f \|^2:= P f^2  $
(again assuming it exists). 
We call $\| f \|_n$ the empirical norm of the function $f$ and
$\| f \|$ its theoretical norm.
We review some results concerning the the uniform (over ${\cal F}$) convergence of $\| \cdot \|_n$ to $\| \cdot \|$.
As example, we consider the case ${\cal X} = \R^p$ (with $p$ possibly large) and ${\cal F}$ is a class of additive
functions $f (x_1 , \ldots , x_p) = \sum_{k=1}^p f_0 (x_k)$ with $f_0$ in a given class of functions ${\cal F}_0$
on $\R$ (Theorem \ref{Sobolev.theorem}). 
We extend the results to uniform convergence of the empirical measure of products of functions. 
The latter will be an important tool for statistical theory for causal inference. As intermediate step we show convergence of the least squares
estimator when the model is wrong.

In Theorem
\ref{square.theorem} we present results from 
\cite{guedon2007subspaces} and \cite{bartlett2012}  and in Theorem \ref{square2.theorem} we compare these with more classical approaches
using e.g. the contraction inequality. The extension
to inner products is given in Theorem \ref{product2.theorem}. 
The latter can be used in statistical applications where functions from different smoothness classes are estimated (for example in an additive model). 

We pay some special attention to the linear case, i.e. the case where ${\cal F}$ is (a subset of) a linear space.
For isotropic distributions the uniform convergence of $\| \cdot \|_n$
to $\| \cdot \|$ over linear functions is well developed. We refer to
\cite{adamczak2011sharp}
 and with sub-Gaussian random vectors to 
 \cite{raskutti2010restricted},  \cite{loh2012} and \cite{rudelson2011reconstruction}. 
 We will not require isotropic distributions but instead consider possibly 
anisotropic but bounded random variables. We present results from \cite{bartlett2012} and 
\cite{rudelson2011reconstruction}
which are based on \cite{guedon2007subspaces} or a similar approach.
Theorems \ref{l1.square.theorem} and \ref{l0.square.theorem} are essentially
in  \cite{bartlett2012} and \cite{rudelson2011reconstruction}.
We compare the bound with a Bernstein type inequality for
random matrices as given in \cite{ahlswede2002strong}.

%The bound we will invoke can be summarized as follows.
%Let ${\cal F} = \{ f_{\beta} ( \cdot) := \sum_{j=1}^p
%\beta_j \psi_j ( \cdot)  : \ \beta \in \R^p \}$ be a linear space in $p$ dimensions. 
%Let  $K_{\psi} := \max_j \| \psi_j \|_{\infty}   $. Furthermore, let $\Lambda_{\rm min}^2$ be
%the smallest eigenvalue of  $\Sigma $, with $\Sigma$ being the population inner product matrix, i.e.
%for $j,k \in \{ 1 , \ldots , p \}$
%$$\Sigma_{j,k} = \sum_{i=1}^n \EE \psi_j (X_i) \psi_k (X_i) / n .$$
%Then 
%$$ \sup_{\beta \in \R^p } \biggl | \| f_{\beta}  \|_n / \| f_{\beta} \| - 1 \biggr | =
%{\mathcal O}_{\PP} \biggl ( \sqrt {K_{\psi}^2 p \log^3 p \log n / (n \Lambda_{\rm min}^2 )} \biggr ) .$$
 
 Uniform convergence of empirical norms and inner products has numerous statistical applications. 
This study is motivated by some questions arising in the structural equations model for causal inference. 
Let us briefly sketch the problem.
Consider having observed an $n \times p$ matrix data matrix $X$ with i.i.d.\ rows.
We assume the structural equations model
$$X_{1,j} = f_{j}^0 \biggl ( {\rm parents} ( X_{1,j} ) \biggr ) + \epsilon_{1,j} , \ j=1 , \ldots , p. $$
Here, ${\rm parents} ( X_{1,j}  ) $ is a subset of $\{ X_{1,k } \}_{k \not= j } $, 
$\epsilon_{1,j}  , \ldots , \epsilon_{1, p} $ are independent Gaussian noise terms,
$\epsilon_{1,j} $ is independent of ${\rm parents} ( X_{1,j} ) $ 
and $f_j^0 $ is the regression of $X_{1,j}$ on its parents ($j=1 , \ldots , p $). 
For a directed acyclic graph (DAG) there exists a permutation
$\pi^0 := ( \pi_1^0 , \ldots , \pi_p^0 ) $ of $\{ 1, \ldots , p \}$ such that for all $j $ the parents of
$X_{1,\pi_j^0}$ are $\{ X_{1, \pi_1^0 }, \ldots ,   X_{1, \pi_{j-1}^0} ,\} $ or a subset thereof, with the
convention that for $j=1$, the parental set is the empty set.
The permutation $\pi_0$ is not unique, and we 
let $\Pi_0$ be the class of permutations with this parental structure.

If for each $j$ the set of parents of $X_{1,j}$ in the DAG were known, the problem is a standard (nonparametric) multiple
regression problem. However, the parental structure, i.e.\ the class $\Pi_0$ is not know and hence
has to be estimated from the data. 
Let  $\Pi$ the class of all $p!$ permutations of $\{  1 , \ldots , p \} $ and 
$\{ {\cal F}_{j} \} _{j=1}^p $ be given classes of regression functions.
Here, ${\cal F}_{j}$ is a collection of functions of $j-1$ variables ($j=1 , \ldots , p$).
We use the short hand notation: for each $i$, $j$ and $\pi$
$$ f_j( X_i , \pi) := f_j (  X_{i, \pi_1} , \ldots , X_{i, \pi_{j-1}} ) , $$
with the above convention for $j=1$, and for each $j$ and $\pi$ 
$$ \| {\bf X}_{\pi_j}  - f_{j }(\pi)  \|_n^2 := \sum_{i=1}^n  (  X_{i,\pi_j } -  f_{j} ( X_i , \pi ) )^2 / n . $$

We consider the estimator 
$$\hat \pi : \in  \arg \min_{\pi \in \Pi} \sum_{j=1}^p \log \biggl ( \| {\bf X} _{\pi_j}  - \hat f_{j} (\pi) \|_n \biggr  )  $$
where, for each $j$, $\hat f_{j }(\pi) $ is the least squares estimator
$$ \hat f_{j} (\pi) := \arg \min_{f_{j} \in {\cal F}_{j }} \| {\bf X}_{\pi_j}  - f_{j} (\pi)\|_n . $$
This estimator is proposed by \cite{PetersB}, where
consistency results, algorithms and simulations are presented.
We further develop the theory using the refined inequalities from
\cite{guedon2007subspaces} and \cite{ahlswede2002strong}. 
We show in Theorem \ref{DAG0.theorem}
that this estimator is consistent under various scenario's: 
$\PP ( \hat \pi \notin \Pi_0)$ converges to zero. An important assumption here is an identifiability assumption: see Condition \ref{identifiable.condition}. This excludes the Gaussian linear structural
equations model where $X_{1,j}$ depends linearly on its parents. We will instead model
each $f_{j }\in {\cal F}_j $ as being an additive non-linear function 
$$f_{j} (x_1 , \ldots , x_{j-1}  )= \sum_{k =1}^{j-1} f_{k,j } ( x_k) ,$$
where each $f_{k,j}$ belongs to a given class ${\cal F}_0 $ of real-valued functions on $\R$.

We consider several cases. The results can be found in Theorem \ref{DAG0.theorem}.
They are a consequence of uniform convergence of empirical norms of a class
of additive functions as given in Theorem \ref{Sobolev.theorem} which may be of independent interest.
Let us summarize the findings here.

 In the first two cases, the class 
${\cal F}_0 $ is assumed to have finite entropy integral for the
supremum norm. We then derive consistency when $p^3 = o(n)$.
Under additional assumptions this is can be relaxed to $p^{3- (1- \alpha)^2} = o(n)$,
where $0< \alpha<1$ is a measure of the ``smoothness " of the class ${\cal F}_0$.

An important special case is where ${\cal F}_0$ is a class of linear functions.
Each $f_{k,j}$ is then a linear combination of functions in a given dictionary $\{ \psi_r \}_{r=1}^N$:
$$f_{k,j} (x_k) = \sum_{r=1}^N \beta_{r,k,j} \psi_r (x_k ) .$$
In other words, the dependence of a variable (index $j$) on one of its parents (index $k$) is then modelled
as a linear combination of certain features (index $r$) of this parent. 
We assume the dictionary to be bounded in supremum norm.

If ${\cal F}_0$ is the signed convex hull of the functions $\{ \psi_r \}_{r=1}^N$
we obtain consistency when $p^2 \log N \log^3 n = o(n) $.
The latter situation covers for example the case where ${\cal F}_0$ is a collection
of functions with total variation bounded by a fixed
constant.

Under certain eigenvalue conditions we find that $pN^2 \log n = o(n)$ also yields consistency.

Finally, if ${\cal F}_0$ can be approximated by linear functions in a space of dimension $N$
with bias of order $N^{-{1 /( 2 \alpha)}}$, then consistency follows from $p^{1 + 4 \alpha} \log n = o(n)$.

The paper \cite{PetersB} shows consistency for the case $p$ fixed (the low-dimensional case).
It also has theoretical results for the high-dimensional case, but for a restricted estimator
where it is assumed that $X_{1,j}$ has only a few parents and a superset of  the parents
${\rm parents} (X_{1,j})$ is known or can be estimated ($j=1 , \ldots , p$). 
This superset then is required to be small.

The paper is organized as follows.
In Sections \ref{empirical.norm.section} and
\ref{empirical.innerproduct.section}
we study a generic class of functions ${\cal F}$
satisfying some $\|  \cdot \|$- and 
$\| \cdot \|_{\infty} $- bounds.
We present the uniform convergence for empirical norms in Section \ref{empirical.norm.section}, with main
example in Subsection \ref{example.section}. 

Section \ref{empirical.innerproduct.section} looks at 
empirical inner products of functions in different ``smoothness" classes.
Subsection \ref{smooth.section} illustrates the results by considering two classes of functions
satisfying different entropy conditions.
In many applications one also needs uniform convergence of
inner products with a sub-Gaussian (instead of bounded) random variable. Therefore we briefly review
this case as well in Subsection \ref{sub-Gaussian.section}.

Section \ref{linear.section} applies  the theory to a class of linear functions and Section \ref{least.squares.section} studies linear regression when the model is wrong. Section \ref{DAG.section} contains the main
application: estimation of the order in a directed acyclic graph.
Section \ref{conclusion.section} concludes.

Section \ref{technical.section} presents the technical tools and Section \ref{proofs.section} contains the proofs.
%We however omit the proofs for Section \ref{linear.section} since here
%straightforward consequences of the preceding sections are presented, inserting an entropy bound given
%in  \cite{rudelson2008sparse} (given in the present paper as Lemma \ref{Rudelson.lemma}). 
Throughout $C_0$, $C_1 $, $C_2$ , $\cdots$ and $c_0$, $c_1$, $c_2$ , $\cdots$ are universal constants, not the same at each appearance.

\section{Bounds for the empirical norm}\label{empirical.norm.section}

\subsection{Entropy and entropy integrals}

For a real-valued function $f$ on ${\cal X}$ we let its supremum norm restricted to the sample be
$$ \| f \|_{n, \infty} := \max_{1 \le i \le n } | f (X_i) | $$ and
we let 
${\cal H} (u, {\cal F} , \| \cdot \|_{n, \infty} ) $ be the entropy of
$({\cal F} , \| \cdot \|_{n , \infty} ) $.
We further define for $z>0$
\begin{equation} \label{J.definition} 
J_{\infty} ^2 (z, {\cal F}) := C_0^2    \inf_{\delta >0} \EE  \biggl [ z  \int_{\delta  }^{1} \sqrt { {\cal H} (uz /2 , {\cal F} , \| \cdot \|_{n, {\infty} } )} du  + \sqrt n \delta z  \biggr ]^2  
\end{equation}
where the constant $C_0$ is taken as in Theorem \ref{Dudley.theorem} (Dudley's Theorem). 
We can without loss of generality assume the  integral exists (replace the entropy by a continuous
upper bound). The subscript $\infty$ here
refers to the fact that we are considering $\ell_{\infty}$-norms. 

We also consider uniform $\ell_2$-entropies, defined as follows.
Let ${\cal A}_n $ be the set of all configurations $A_n$ of $n$ (possibly non-distinct) points within the support of $P$.
For $A_n \in {\cal A}_n$ and $f$ a real-valued function on ${\cal X}$  we let
$$ \| f \|_{A_n}^2 := \sum_{x \in A_n }  f^2(x) / n . $$
Note that $\| f \|_{n } = \| f \|_{{\bf X} } $ where ${\bf X}$ is the random sample ${\bf X} := \{X_1 , \ldots , X_n \} $.
For a class ${\cal F}$ of functions on ${\cal X}$, we let
$${\cal H} ( \cdot , {\cal F}) := \sup_{A_n \in {\cal A}_n } {\cal H} ( \cdot , {\cal F} , \| \cdot \|_{A_n } )  $$
and 
\begin{equation}\label{calJ.definition}
{\cal J}_0( z , {\cal F} ) := C_0 z \int_{0}^1 \sqrt {{\cal H} ( uz /2, {\cal F} )} du 
  , \ z > 0 . 
  \end{equation}
  The calligraphic symbol ${\cal J}$ indicates that instead of random entropies we consider the maximum entropy 
  over all possible configurations of (at most) $n$ points. Apart from this and from  considering $\ell_2$-entropy
  instead of $\ell_{\infty}$-entropy we now
moreover implicitly assume that the entropy integral converges and use ${\cal J}_0$ with subscript 0
to indicate this. The reason for taking 0 as lower-integrant is that $v \mapsto
{\cal J}_0 ( \sqrt {v}, {\cal F}) $ is a concave function. We will see this to be useful in Theorem \ref{square2.theorem} in view
of Jensen's inequality.

Finally, for $A_n \in {\cal A}_n$ and $f$ a real-valued function on ${\cal X}$  we let
$$ \| f \|_{A_n, {\infty}} := \max_{x\in A_n } | f(x) | . $$
Note that $\| f \|_{n, {\infty} } = \| f \|_{{\bf X} , \infty} $ where ${\bf X}$ is the sample ${\bf X} := \{X_1 , \ldots , X_n \} $.
For a class ${\cal F}$ of functions on ${\cal X}$ we set
$${\cal H}_{{\infty}} ( \cdot , {\cal F}) := \sup_{A_n \in {\cal A}_n } {\cal H} ( \cdot , {\cal F} , \| \cdot \|_{A_n , 
\infty} ) . $$
We furthermore define for $z > 0$ 
\begin{equation}\label{Jinfty.definition}
{\cal J}_{\infty} ( z , {\cal F} ) := C_0  \inf_{\delta >0 } \biggl [z  \int_{\delta  / 4}^1 \sqrt {{\cal H}_{ {\infty} } ( uz/2, {\cal F} )} du +
\sqrt {n} \delta z  \biggr ]  . 
\end{equation}
By the definition of $J_{\infty} $ (see (\ref{J.definition}))
$J_{\infty} (z, {\cal F}) \le {\cal J}_{\infty} ( z , {\cal F})$.
We use the calligraphic symbol ${\cal J}_{\infty}$ with subscript $\infty$ here to indicate that the maximal
$\ell_{\infty}$-entropy over all possible configurations of (at most)  $n$ points is used.
  
  \subsection{Bounds using $\ell_{\infty}$-norms}\label{linfty.section}

The following theorem follows from \cite{guedon2007subspaces}. Recall the definition (\ref{J.definition})
of $J_{\infty}$. 
\begin{theorem} \label{square.theorem} Let
$$ R:=  \sup_{f \in {\cal F} } \| f   \|  , \ K:= \sup_{f \in {\cal F}} \| f \|_{\infty} . $$
Then
$$ \EE \biggl (\sup_{f \in {\cal F} } \biggl |  \| f \|_n^2 - \| f \|^2 \biggr | \biggr ) \le 2 J_{\infty} (K, {\cal F}) R / \sqrt n+
4 J_{\infty}^2 (K, {\cal F })/n . $$
Moreover, for all $t >0$, with probability at least $1- \exp[-t]$, 
$$\sup_{f \in {\cal F}} \biggl | \| f \|_n^2 - \| f \|^2 \biggr | / C_1 \le
{2 R J_{\infty} (K, {\cal F} ) +RK \sqrt t   \over  \sqrt n  }    + { 4  J_{\infty}^2 (K, {\cal F})  + K^2 t \over n}     $$
where the constant $C_1$ is as in Theorem \ref{bounded.theorem} (a deviation inequality).  \\
As by-product of the proof, we find
$$\sqrt  {\EE \hat R^2 } \le R + 2 J_{\infty} (K, {\cal F} )/\sqrt n  . $$

\end{theorem}

Actually, in \cite{guedon2007subspaces} the entropy integral related quantity $J_{\infty}$ is replaced
by a more general quantity coming from generic chaining. 

\subsection{Bounds using $\ell_2$-norms}\label{l2.section}

In Theorem \ref{square2.theorem} below, we reverse the role of $R$ and $K$ as compared to
Theorem \ref{square.theorem}.
The result is well-known, it follows from contraction inequality (\cite{Ledoux:91})
or from a direct argument. See also \cite{gine2006concentration}.
Recall the definition (\ref{calJ.definition}) of ${\cal J}_0$. 

\begin{theorem} \label{square2.theorem} Let
$$ R:=  \sup_{f \in {\cal F} } \| f   \|  , \ K:= \sup_{f \in {\cal F}} \| f \|_{\infty} . $$
Let for $z>0$, $G^{-1} (z^2) := {\cal J}_0 (z, {\cal F})  $ and let $H$ be the convex conjugate of $G$.
Assume that $R^2 \ge H( 4 K / \sqrt n ) $. 
Then
$$ \EE \biggl (\sup_{f \in {\cal F} } \biggl |  \| f \|_n^2 - \| f \|^2 \biggr | \biggr ) \le { 2K {\cal J}_0  (2R, {\cal F})  \over  \sqrt n} 
 . $$
Moreover, for $R^2 \ge H( 4 K / \sqrt n ) $ and all $t >0$
$$\PP \biggl (\sup_{f \in {\cal F}} \biggl | \| f \|_n^2 - \| f \|^2 \biggr | / C_1 \ge { 2 K {\cal J}_0  (2R
, {\cal F}) + KR \sqrt t  \over \sqrt n }+
  {  K^2 t \over n}  \biggr )  \le \exp[-t]  $$
where the constant $C_1$ is as in Theorem \ref{bounded.theorem}. \\
As by-product of the proof, we find
$$\EE \hat R^2 \le 4R^2  . $$

\end{theorem}

%Theorem \ref{square2.theorem} is proved by direct arguments, but can also be obtained using
%the contraction inequality of ..., We give here an intermediate result of this - more classical - approach,
%to show in Example ... that it can actually lead to better bounds than Theorem \ref{square.theorem}
%when $J_{\infty} (K, {\cal F} )/ K$ is very large.
%
%\begin{lemma} \label{Rademacher.lemma} Let $\{ \epsilon_i \}_{i=1}^n $ be a Rademacher sequence, independent of
%$X_1, \ldots , X_n $ and define
%$$P_n^{\epsilon} f := {1 \over n} \sum_{i=1}^n \epsilon_i f(X_i) . $$
%Then for $R^2 \ge H( 4K/ \sqrt n) $ and all $t \ge 0 $
%$$ \PP \biggl (\sup_{f \in {\cal F}} \biggl |P_n^{\epsilon} f  \biggr | / C_1 \ge {  {\cal J}_0  (2R
%, {\cal F}) + R \sqrt t  \over \sqrt n }+
%  {  K t \over n}  \biggr )  \le \exp[-t]  $$
%where the constant $C_1$ is as in Theorem \ref{bounded.theorem}.
%\end{lemma} 
%
%\begin{corollary} \label{Rademacher.corollary}Suppose that for all $z >0$, ${\cal J}_0 (z, {\cal F} ) \le c_0 z^{1- \alpha}$ for some
%fixed $0 < \alpha < 1$. Assume moreover that $R=1$ and that $\| f \|_{\infty} / \| f \|^{1- \alpha} \le c_1 $
%for all $f \in {\cal F} $. Then, by the peeling device ...
%$$\PP \biggl (\sup_{f \in {\cal F}} \biggl |{ P_n^{\epsilon} f  \over
% \| f \|^{1- \alpha} + n^{-{1- \alpha  \over 2(1+\alpha) }} }  \biggr | / c_2 \ge  {1+ \sqrt t \over \sqrt n} + {t \over n}  \biggr )  \le \exp[-t]  $$

\subsection{The scaling phenomenon} \label{scaling.section}
As said, the essential difference between Theorems \ref{square.theorem} and \ref{square2.theorem}
is that the roles of $K$ and $R$ are reversed, instead of $R J_{\infty} (K,{\cal F})$ we are
dealing with $K {\cal J}_0 ( R, {\cal F} ) $. In some situations $J_{\infty} (K, {\cal F} )/ K$ behaves as a constant
whereas ${\cal J}_0 (R, {\cal F})/R $ decreases in $R$. Let us illustrate this here. 
Let ${\cal F}_1$ be a class of functions, uniformly $\| \cdot \|_{\infty}$-bounded  by 1, and
consider for $R \le 1 $ the localized class
$${\cal F} (R) := \{ f \in {\cal F}_1: \ \| f \| \le R \} . $$
Suppose that 
$$\sup_{f \in {\cal F}_1 (R) } \| f \|_{\infty} \asymp 1 ,  \ 0 < R \le 1 $$
and for some $0< \alpha < 1$
$$J_{\infty} (z, {\cal F}_1 (R) )  \asymp 
{\cal J}_0 (z , {\cal F}_1 (R)  )  \asymp z^{1- \alpha}, \ z > 0 ,  \ 0 < R \le 1. $$
These assumptions say that the local class ${\cal F}(R)$ behaves like to global class
${\cal F}_1$ as far as supremum norm and entropy are concerned. 
Then, taking $K\asymp 1$,
$$R J_{\infty} (K, {\cal F} (R)) \asymp  R, \ K {\cal J}_0  (R , {\cal F} (R) ) \asymp R^{1- \alpha} . $$
Thus, by using Theorem \ref{square.theorem} instead of Theorem \ref{square2.theorem} we win
a factor $R^{\alpha}$.

Otherwise put, let ${\cal F}_K := \{ Kf : \ f \in {\cal F}_1 \} $ for some $K \ge 1 $ and
$${\cal F}_K (1) := \{ f \in {\cal F}_K:\ \| f \| \le 1 \} . $$
Then, taking $R=1$,
$$R J_{\infty} (K, {\cal F}_K (1) ) \asymp K ,\ K {\cal J}_0 ( R , {\cal F}_K (1)) \asymp K^{1+\alpha} .  $$
So by using Theorem \ref{square.theorem} instead of Theorem \ref{square2.theorem} we get rid of
a factor $K^{\alpha}$.

In fact, we find a scaling phenomenon in Theorem \ref{square.theorem}:
whereas for general deviation inequalities the term involving the expectation of the supremum
of the empirical process dominates the deviation term, in the current situation they are of the same order.

Also more generally Theorem \ref{square.theorem} gives better results than Theorem
\ref{square2.theorem}. As we will see, in the particular case where ${\cal F}$ is the signed convex hull of $p$
given functions, uniform convergence follows from Theorem \ref{square.theorem} for $p$
of small order $n$ (up to log-factors) (see Theorem \ref{l1.square.theorem}), whereas Theorem \ref{square2.theorem} needs $p$ to be
of small order $\sqrt n$ (up to log-factors).

 \subsection{Example: additive functions} \label{example.section} 
 Let ${\cal F}_0$ be a class of real-valued functions defined on the real line. 
 Let further ${\cal X}:= \R^p$  where $ p \le n$ and let
 $${\cal F} := \biggl \{ f (x_1 , \ldots , x_p ) = \sum_{k=1}^p f_{k}  ( x_k) : \ f_k\in {\cal F}_0 \ \forall \ k \biggr  \} . $$
 
 We will sometimes require the following incoherence condition: for a constant $c_1$
 and for all $f_0 \in {\cal F}_0$ and $f_{0,k}(x_1 , \ldots , x_p ) := f_0 (x_k)$, $k=1 , \ldots , p $, 
 \begin{equation}\label{incoherence.equation}
 \sum_{k=1}^p \| f_{0,k} \|^2 \le c_1 \| \sum_{k=1}^p f_{0,k} \|^2  .
 \end{equation}
 
 In the following theorem one may think of
 ${\cal F}_0$ being for a given $m \in \Nat$ the Sobolev class.
 \begin{equation}\label{Sobolev.equation}
 {\cal F}_0 = \biggl \{ f_0 : \ \int | f_0^{(m)} (v ) |^2 dv \le 1 \biggr \} .
 \end{equation}  
The constant $\alpha$ is then $\alpha = 1/ (2m)$ and the choice $N \asymp 
n^{ \alpha \over 1+ \alpha}$ corresponds to taking a piecewise polynomial
approximation with $\asymp n^{1 \over 2m+1} $ pieces (i.e. the bandwidth
of the usual order $n^{-{m \over 2m+1}}$).  The bound (\ref{supnorm.equation}) is shown
for this case in \cite{Agmon:65} under  the condition that the one-dimensional marginal densities
of the $X_i \in \R^p$ stay away from zero (see also Lemma \ref{Agmon.lemma} below).

We define
$${\bf Z} ({\cal F} (1)) := \sup_{f \in {\cal F}, \ \| f \| \le 1 } \biggl | \| f \|_n^2 - \| f \|^2 \biggr | . $$
 
 \begin{theorem} \label{Sobolev.theorem} $ $\\
 {\bf Case 1.} 
  Assume that for a fixed $0 < \alpha < 1$, 
  \begin{equation}\label{entropyalpha.equation}
 \int_0^z \sqrt { H( u , {\cal F}_0 , \| \cdot \|_{\infty} ) } du   \asymp z^{1- \alpha},  \ z > 0 .
 \end{equation}
 Then ${\bf Z}^2 ({\cal F} (1)) = {\mathcal O}_{\PP} ( {p^3 / n }) $.\\
 {\bf Case 2.} Assume in addition to the condition of Case 1 that the incoherence condition
 (\ref{incoherence.equation}) holds true for some constant $c_1 = {\mathcal O} (1)$
 and that for some constant $c_2 = {\mathcal O} (1)$
 and for all $f_0 \in {\cal F}_0$, all $j$, and for $f_{0,k} (x_1 , \ldots , x_p ) = f_0 (x_k)$
 \begin{equation}\label{supnorm.equation}
 \| f_{0,k} \|_{\infty} \le c_2 \| f_{0,k} \|^{1- \alpha} .
 \end{equation}
 Then ${\bf Z}^2 ({\cal F} (1)) = {\mathcal O}_{\PP} ( {p^{3- ( 1 - \alpha )^2 } / n} ) $,\\
 {\bf Case 3.} Suppose that ${\cal F}_0$ is the signed convex hull of given functions $\{ \psi_r \}_{r=1}^N$,
 $N \le n$,  in
 particular
 \begin{equation}\label{convexhull.equation}
 {\cal F}_0 = \biggl \{ f_0 = \sum_{r=1}^N \beta_r \psi_r ( \cdot ) : \ \sum_{r=1}^N | \beta_r | \le 1 \biggr \} , 
 \end{equation}
 where $\{ \psi_r \}$ is a given dictionary satisfying $\max_r \| \psi_r \|_{\infty} = {\mathcal O} (1)$.
 Then ${\bf Z}^2 ({\cal F}(1) ) = {\mathcal O}_{\PP} (  {p^2 \log^4 n  / n} ) $.\\
 {\bf Case 4.} Suppose that for some $N \le n $
 $${\cal F}_0 = \biggl \{ f_0 (\cdot) = \sum_{r=1}^N \beta_r \psi_r ( \cdot ): \  ( \beta_1 , \ldots , \beta_N ) \in \R^N 
 \biggr \}  , $$
 where $\{ \psi_r \}$ is a given dictionary satisfying $\max_r \| \psi_r \|_{\infty} = {\mathcal O} (1)$.
 Assume that the incoherence condition (\ref{incoherence.equation}) is met for some constant $c_1 = {\mathcal O} (1)$.
 Assume moreover that for  a constant $c_0 = {\mathcal O} (1)$, all $\beta \in \R^N$, and for all $k$ and for
 $f_{\beta,k} (x_1 , \ldots , x_p) := \sum_{r=1}^N \beta_r \psi_r (x_k) $, 
 \begin{equation}\label{eigenvalue.equation}
 \| \beta \|_2^2 \le c_0 N \| f_{\beta , k} \|^2 . 
 \end{equation}
 Then ${\bf Z}^2 ({\cal F}(1) ) = {\mathcal O}_{\PP} (  {p N^2   \log n/ n} ) $.
 When one chooses $N \asymp n^{\alpha \over  1 + \alpha} $
 this reads
 ${\bf Z}^2 ({\cal F} (1)) = {\mathcal O}_{\PP} ( {p n^{\alpha \over 1+ \alpha}    \log n / n^{1 \over 1+ \alpha} } ) $.
  \\
 {\bf Case 5.} Consider  a dictionary $\{ \psi_r \}_{r=1}^{\infty}$ with $\sup_r \| \psi_r \|_{\infty}
 = {\mathcal O} (1)$.
 Suppose that for a constant $c_0 = {\mathcal O} (1)$ and any $f_0 \in {\cal F}_0$
 and any $N \in \Nat$ there exists a $\beta \in \R^N$ such that
 $$\biggl \|  f_0 - \sum_{r=1}^N \beta_r \psi_r \biggr \|_{\infty} \le c_0 N^{-{ 1  \over 2 \alpha}}  . $$
 Moreover, assume the incoherence condition
 (\ref{incoherence.equation}) with $c_1 = {\mathcal O}(1)$and that for all $N \in \Nat$, all $\beta \in \R^N$, all $k$ and for
 $f_{\beta,k} (x_1 , \ldots , x_p) := \sum_{r=1}^N \beta_r \psi_r (x_k) $, 
 $$ \| \beta \|_2^2 \le c_0 N \| f_{\beta , k} \|^2 . $$
Then  ${\bf Z}^2 ({\cal F} (1)) = {\mathcal O}_{\PP} ({p^{1+4 \alpha \over 1 + 2 \alpha } (\log n  / n)^{1 \over 1+ 2 \alpha } }) $.\\
  \end{theorem}
  
  \begin{remark} \label{local.remark} If Condition (\ref{incoherence.equation}) holds, one may in fact replace
  ${\cal F}_0 $ in (\ref{entropyalpha.equation}) of  Case 1 by the local class ${ \cal F}_0 \cap \{  \min_k \| f_0\|_k^2 \le c_1 \} $
  where $\| f_0 \|_k^2= \EE \sum_{i=1}^n f_0^2 (X_k)  /n$:
 $$ \int_0^z \sqrt { H( u , {\cal F}_0 \cap \{ \min_k  \| f_0\|_k^2 \le c_1 \}  , \| \cdot \|_{\infty} ) } du   \asymp z^{1- \alpha},  \ z > 0  .$$ 
 The same is true for Case 2, where Condition (\ref{incoherence.equation}) is indeed assumed.
 In Case 3, assuming (\ref{incoherence.equation}) one may replace condition (\ref{convexhull.equation}) by the local version
 $${\cal F}_0 \cap \{ \min_k \| f_0 \|_k^2 \le c_2 \} \subset
 \biggl \{ f_0 = \sum_{r=1}^N \beta_r \psi_r ( \cdot ) : \ \sum_{r=1}^N | \beta_r | \le 1 \biggr \} .$$

 \end{remark} 
  
  To complete the picture we show in the next lemma that condition (\ref{supnorm.equation}) is natural
  in the context of Case 5 (although we do not use it there).
  
  \begin{lemma}\label{Agmon.lemma} Assume the conditions of Case 5 in Theorem \ref{Sobolev.theorem} and 
  that 
  for some constant
 $c_3 = {\mathcal O} (1)$, for all $r$ and all $s > c_3$,
 \begin{equation}\label{nonoverlap.equation}
 \psi_r \psi_{r+s} =0.
 \end{equation} 
i.e., that as soon as $s > c_0$, $\psi_r $ and $\psi_{r+s}$ do not overlap.Then (\ref{supnorm.equation}) holds for some
constant $c_2 = {\mathcal O}(1)$.
  \end{lemma}

%whereas Theorem \ref{square2.theorem} gives
%$$\sup_{f \in {\cal F} , \ \| f \| \le 1 } 
% \biggl | \| f \|_n ^2 - \| f \|^2 \biggr | = {\mathcal O}_{\PP} (
% p^{2 + 3 \alpha \over 2   } )/  \sqrt n . $$
%The first bound is always better than the last one.
In Case 2, the bound found in \cite{meier2009high} is ${\bf Z}( {\cal F}(1))= {\cal O}_{\PP}(  p^{2(1+\alpha)} /  n ) $.
Note that in Case 5, we have ${\bf Z} ({\cal F} (1)) = o_{\PP} (1)$ whenever
$p^{1+ 4 \alpha} / n = o (1)$. 
The conditions on $p$ can possible be weakened (possibly by replacing entropy
bounds by Gaussian means) but this is an open problem.
It is not clear to us whether the bounds presented in Theorem \ref{Sobolev.theorem}
are sharp.
% However, by the contraction inequality, we see that
% $$ \sup_{f \in {\cal F} , \| f \| \le 1 } 
% \biggl | \| f \|_n ^2 - \| f \|^2 \biggr | = p^{1+ \alpha \over 2}  \sup_{f \in {\cal F} , \| f \| \le 1 } 
% \biggl | P_n^{\epsilon} f  \biggr | {\mathcal O}_{\PP} (1) = p^{1+ \alpha } . $$

Case 1 and 2 of Theorem \ref{Sobolev.theorem} follow from Theorem
\ref{square.theorem} by straightforward entropy bounds.
Case 3 is based on  a result from
\cite{rudelson2011reconstruction} cited here as Theorem \ref{l1.square.theorem}.
Case 4 is based on the general matrix version of Bernstein's inequality of \cite{ahlswede2002strong} cited
here as Theorem \ref{l3.square.theorem}.  Case 5 follows from Case 4 using a trade-off argument for
the choice of $N$ (the value $N= n^{\alpha / (1+ \alpha)}$ suggested in Case 4 may not give the optimal trade-off).
 The details are in Section \ref{proofs.section}.

\section{Empirical inner products}\label{empirical.innerproduct.section}

%\begin{lemma} Let
%$$J(z, {\cal F}) := z \int_{0}^{1} \sqrt {H(uz/2 , {\cal F} , \| \cdot \|_{\infty} ) du }:= G^{-1} (z^2) . $$
%Assume that $G$ is strictly convex. Let $H$ be its convex conjugate. Then
%for any random variable $Z$
%$$\EE J(Z, {\cal F}) \le G^{-1} \EE(Z^2) .$$
%
%\end{lemma}
%

Consider products $fg$  of functions $f$ and $g$ with $f$ in some class ${\cal F}$ and $g $ in some class
${\cal G}$. Note that one can derive results for products via squares:
$$fg = (f+g)^2 /2 - (f^2 + g^2 )/ 2 . $$
If ${\cal F}$ and ${\cal G}$ have the same $\| \cdot \| $-diameter $R$ and the same
$\| \cdot \|_{\infty}$-diameter $K$ it is easy to see that without loss of generality
we may assume that ${\cal F} = {\cal G}$ (replace ${\cal F}$ and ${\cal G}$ by ${\cal F} \cup {\cal G} $).
%This case is studied in Subsection \ref{same.section}.
However, if $f$ and $g$ are in different classes it may be more appropriate to analyze the products directly.
This case with ${\cal F}$ and ${\cal G}$ having different radii is studied here.

We only present the results using $\ell_{\infty}$-norms. Again, one may reverse the roles
of $\| \cdot \|_{\infty}$-radii and $\| \cdot \|$-radii, getting other versions for the bounds.
The best bound may depend on the situation at hand. 

%In fact, from Remark
%\ref{singleton.remark} we see that there is no longer a clear reversion of roles.

%\subsection{Inner products of functions from the same class}\label{same.section}
%
%We recall the definition of $J_{\infty} (z, {\cal F})$ given in (\ref{J.definition}).
%
%\begin{theorem} \label{product.theorem} Let
%$$ R:= \sup_{f \in {\cal F}} \| f \| , \ K:= \sup_{f \in {\cal F}}  \| f \|_{\infty} . $$
%For all $t>0$, with probability at east $1- \exp[-t]$
%$$\PP \biggl (\sup_{f , g \in {\cal F}} \biggl | (P_n - P) fg   \biggr | / C_1 \ge
%{4   R J_{\infty} (K, {\cal F} )  + RK \sqrt t  \over  \sqrt n  }  + {8    J_{\infty} ^2 (K, {\cal F}) + K^2 t  \over n} 
%\biggr ) $$
%where the constant $C_1$ is as in Theorem \ref{bounded.theorem}. 
%\end{theorem}
%
%Note that indeed Theorem \ref{product.theorem} gives up to a factor 2 exactly the same bound as
%Theorem \ref{square.theorem}.
%
%
%
%
%
\subsection{Inner products of functions from different classes}\label{different.section} 
%
%Consider now two classes of functions ${\cal F}$ and ${\cal G}$. 
Let 
$$ R_1 := \sup_{f \in {\cal F}} \| f \|   , \ K_1 := \sup_{f \in {\cal F} }\| f \|_{\infty} . $$
and 
$$R_2 := \sup_{g\in {\cal G}} \| g \|  , \ K_2 := \sup_{g \in {\cal G} }\| g \|_{\infty} .$$

\begin{theorem} \label{product2.theorem} Suppose that $R_1 K_2 \le R_2 K_1 $. Consider values of $t\ge 4$ and $n$ such that
\begin{equation} \label{R1.equation}
\biggl ( {2 R_1  {\cal J }_{\infty} (K_1, {\cal F} ) + R_1 K_1 \sqrt t   \over  \sqrt n  }  + { 4  {\cal J}_{\infty}^2 (K_1, {\cal F}) + K_1^2 t  \over  n}   \biggr ) \le {R_1^2  \over  C_1} 
  \end{equation}
  and
  \begin{equation} \label{R2.equation}
 \biggl ( {2 R_2  {\cal J}_{\infty}  (K_2,  {\cal G} )  + R_2 K_2 \sqrt t  \over  \sqrt n  }  + { 4  {\cal J}_{\infty} ^2 (K_2 , {\cal G}) + K_2^2 t  \over  n} 
   \biggr ) \le {R_2^2  \over  C_1} .
   \end{equation}
Then with probability at least $1 - 12 \exp [-t]$
$${1 \over 8 C_1} \sup_{f \in {\cal F} , \ g \in {\cal G}} \biggl | (P_n- P) fg\biggr | \le 
 { R_1 {\cal J}_{\infty} ( K_2 , {\cal G} )  +  R_2 {\cal J}_{\infty}  ( R_1 K_2  / R_2, {\cal F} ) + R_1 K_2 \sqrt t \over  \sqrt n} 
  $$ $$  +  { K_1 K_2  t \over n } 
 . $$
\end{theorem}

\begin{remark} Theorem \ref{product2.theorem} can be refined using generic chaining type of quantities
instead of entropies. We have omitted this to avoid digressions.
\end{remark}

\begin{remark} \label{singleton.remark} Consider the special case where ${\cal G} = \{ g_0 \}$ is a singleton.
Assume that $\| g_0 \|_{\infty} =K_0$. Take $R_2=K_2=K_0$ in Theorem \ref{product2.theorem},
and write $R_1:=K$ and $K_1 := K$. For a singleton ${\cal G}$, the term ${\cal J}_{\infty} (K_2 , {\cal G})$
can be omitted. We then get from Theorem \ref{product2.theorem}: for
$t \ge 4 $ and
$$ \biggl ( {2  R {\cal J }_{\infty} (K, {\cal F} )  +RK \sqrt t  \over   \sqrt n  }  + { 4  {\cal J}_{\infty}^2 (K, {\cal F}) 
+ K^2 t \over  n}   \biggr ) \le {R^2  \over  C_1} ,$$
  it holds that
  $${1 \over 8 C_1} \sup_{f \in {\cal F}} \biggl | (P_n- P) fg_0\biggr | \le 
  {  K_0 {\cal J}_{\infty}  ( R , {\cal F} ) + K_0 R \sqrt t \over  \sqrt n} 
  +  { K K_0  t \over n } 
 $$
  with probability at least $1-8 \exp[-t]$.
 We will see a similar result in Theorem \ref{subGaussianproduct.theorem}, where $g_0$ is not bounded but sub-Gaussian.

  \end{remark}

\subsection{Empirical inner products for smooth functions}\label{smooth.section}

Let us suppose that
$${\cal J}_{\infty} (z,  {\cal F} ) \asymp z^{1- \alpha} , \ {\cal J}_{\infty} (z, {\cal G}) \asymp z^{1- \beta} , $$
where $\beta > \alpha$. For example, one may think of Sobolev classes as
was indicated in Subsection \ref{example.section}, or more locally adaptive
cases such as  ${\cal F} = \{f: [0,1] \rightarrow [0,1] ,\  \int | f^{\prime \prime } (x)  | dx \le 1 \} $
and ${\cal G} \subset \{ g : [0,1] \rightarrow [0,1] :\ \int | g^{\prime} (x) | dx \le 1 \} $. 
Then ${ \cal J}_{\infty} (z,  {\cal F} ) \asymp z^{3/4}$ ($\alpha =1/4$) and ${ \cal J}_{\infty} (z,  {\cal G} ) 
\asymp z^{1/2} \sqrt {\log n} $ ($\beta = 1/2$). The $\log n$-term plays a moderate role and we neglect such
details in the following general line of reasoning.

The fact that $\beta > \alpha$ expresses that ${\cal F}$ is smoother (less rich) than ${\cal G}$.
Having an additive model in mind (the response $Y_i $ is an additive function plus noise
$Y_i= f^0(X_{i,1}) + g^0(X_{i,2}) + \varepsilon_i$, $i=1  , \ldots , n$)
one may expect to be able to estimate a function $f^0 \in {\cal F}$ with squared rate
$R_1^2 := n^{-{1 / (1+ \alpha) } } $ and a function $g^0 \in {\cal G} $ with (slower) squared rate
$R_2^2 := n^{-{1 /  (1+ \beta) } } $.  Let us simplify the situation by assuming that
$X_{i,1}$ and $X_{i,2}$ are independent (the dependent case is detailed in \cite{vdGMammen13}). 
Also assume that the functions in ${\cal F}$ and ${\cal G}$ are already centred. 
We now want to show that $P_n fg $ is small, namely negligible as compare to $R_1^2$.
Indeed, inserting Theorem \ref{product2.theorem} (note that (\ref{R1.equation}) and
(\ref{R2.equation}) are true for $t$ fixed and $n$ sufficiently large), we get with probability at least
$1- 12\exp[-t]$ 
$$ \sup_{f \in {\cal F} , \| f \| \le R_1  , \ g \in {\cal G} , \ \| g \| \le R_2 }{ | P_n  fg | /c_1 \over R_1^2 }   \le \biggl (  {1 \over
\sqrt n R_1} + R_2^{\alpha}  +   \sqrt { t \over  n R_1^2 } +{ t \over n R_1^2 } \biggr ) . $$
For fixed $t$ the right hand side of the above inequality is $o(1)$. 

Actually, \cite{vdGMammen13} first proof  the global (slow) rate $R=R_2$. Suppose that
that now $f/K_1 \in {\cal F} $ where $K_1 = R / \lambda $ with $\lambda \asymp
n^{-1 / (1+ \alpha) }$. 
Again 
(\ref{R1.equation}) and
(\ref{R2.equation}) are true for $t$ fixed and $n$ sufficiently large for $R_1^2=R_2^2 = R^2= n^{-1/(1+\alpha)}$, 
$K_1= R / \lambda$ and $K_2 =1$. 
We find as similar result as above: with probability at least $1- 12\exp[-t]$
$$ \sup_{f/ K_1  \in {\cal F} , \ \| f \| \le R  , \ g \in {\cal G} , \ \| g \| \le R }{ | P_n  fg | /c_1 \over R^2 }  \le {1 \over \sqrt n R} + {  R^{\alpha}  } + \sqrt { t \over n R^2} +
 { t \over n \lambda } .$$

Related is the paper \cite{MullervdGeer13} where the additive model is studied with
$f^0$ a high-dimensional linear function. Again, it can be shown that $f^0$ can be estimated
with a fast oracle rate, faster than the rate of estimation of the unknown function $g^0$.

\subsection{Products with a sub-Gaussian random variable}\label{sub-Gaussian.section}
Consider now real valued random variables $Y_i$, $i=1, \ldots , n $.
We let $P_n$ be the empirical measure based on $\{ X_i , Y_i \}_{i=1}^n $: 
for a real-valued function $f$ on ${\cal X}$
$$P_n {\bf Y} f:= \sum_{i=1}^n Y_i f(X_i) / n .$$
We write $P{\bf Y} f := \EE P_n {\bf Y} f $. 
We study the supremum of the absolute value of the product process
$( P_n - P) {\bf Y} f , \ f \in {\cal F} $. 

\begin{definition} For $Z \in \R$ and $\Psi(z) := \exp [ |z|^k ] $, $k=1 , 2$, we define the
Orlicz norm
$$\| Z \|_{\Psi_k} := \inf \{ L>0 : \EE \Psi_k (Z/L) -1 < 1 \} ,$$
whenever it exists. If $\| Z \|_{\Psi_1}$ exists, we call $Z$ sub-exponential,
and if $\| Z \|_{\Psi_2}$ exists we call $Z$ sub-Gaussian.
\end{definition} 

\begin{definition} \label{conditional-subGaussian.definition} We say that ${\bf Y}:= \{ Y_1 , \ldots , Y_n \}$ is uniformly sub-Gaussian
with constant $K_0$ if 
 $$ \max_{1 \le i \le n } \| Y_i \|_{\Psi_2} \le  K_0 < \infty . $$
\end{definition}

The result below is about products of functions, where the class ${\cal G}$ consists of the single
sub-Gaussian function ${\bf Y}$. 

We recall the definition (\ref{Jinfty.definition}) of ${\cal J}_{\infty} $.

\begin{theorem} \label{subGaussianproduct.theorem} 
Let 
$$
 \sup_{f \in {\cal F}} \| f \| \le R , \ K:= \sup_{f \in {\cal F}}  \| f \|_{\infty} . $$
 Suppose ${\bf Y}$ is uniformly sub-Gaussian with constant $K_0$.
 Consider values of $t$ and $n$ such that
$$ \sqrt {2t \over n} + {t \over n } \le 1 . $$
For these values
$$\PP \biggl (\sup_{f \in {\cal F}}  | (P_n - P) {\bf Y} f  |/ C_1 \ge
 {   2 {\cal J}_0 (K K_0 , {\cal F} )  + K K_0 \sqrt t \over  \sqrt n  }   
  \biggr )  \le 8 \exp[-t] $$
where the constant $C_1$ is  as in Theorem \ref{subGaussian.theorem}. 
\end{theorem}

%The implicit assumption in the definition (\ref{calJ.definition}) of ${\cal J}_0$ that the entropy integral converges
%is not crucial here: we may also include a lower integrant. The argument for this is based on the 
%Cauchy-Schwarz inequality $| P_n Y f | \le \| Y \|_n  \| f \|_n $.  We omit the details but only
%note that if $\| Y_i \|_{\Psi_2 } \le K_0$ for all $i$, then the squared random variable
%$Y_i^2$ is sub-exponential for all $i$. So then by Bernstein's inequality, for all $t>0$,
%\begin{equation} \label{subexponential.equation}
%\PP \biggl ( | \| Y \|_n^2 - \| Y \|^2  | / C_1 \ge K_0^2 \sqrt {t/n} + K_0^2 t/n \biggr ) \le \exp[-t] . 
%\end{equation}
%

\section{Application to a class of linear functions}\label{linear.section}
Suppose ${\cal X} = \R^p$. We let $X_i$ be a row vector in $\R^p$, $i=1 , \ldots , n $.
For a column vector $\beta \in \R^p $ we define $f_{\beta}  (X_i) := X_i \beta $. 
We assume in that for some constant $K_X$
$$\max_{i, j} | X_{i,j} | \le K_X . $$

The following lemma is Lemma 3.7 in \cite{rudelson2008sparse}. We  inserted an explicit
constant.
\begin{lemma} \label{Rudelson.lemma} We have
$${\cal H} (u, \{f_{\beta}: \ \| \beta \|_1 \le 1 \} , \| \cdot \|_{n, \infty } ) \le \biggl (  1+{8 \log (2p) \log(2n) K_X^2  \over u^2 } 
 \biggr ) , \ u > 0 . $$

\end{lemma}

%\begin{lemma}  \label{l1-entropy.integral.lemma} We have 
%$$J( \{f_{\beta} : \  \| \beta \|_1 \le M \} ) \le 2C_0 \sqrt {\log (2n)\log (2p) } MK_X \log n +C_0 MK$$
%where the constants $C_0$ is as in Theorem \ref{Dudley.theorem}.
%
%\end{lemma}
%
%{\bf Proof.}
%If $\| \beta \|_1 \le M $ we know that 
%$$\| f_{\beta} \|_{\infty} \le \|\beta \|_1 
%\max_{1 \le i \le n } \max_{1 \le j \le p } |X_{i,j} | \le MK_X . $$
%Fixing $\delta $ at $\delta = 1/\sqrt n $ in (\ref{J.definition}), we find by Lemma \ref{Rudelson.lemma}
%$$J({\cal F}) \le C_0 MK_X
%\int_{1/  \sqrt {n}}^{1} \sqrt { {\cal H} (u/2 , \{f_{\beta} : \  \| \beta \|_1 \le M \} , \| \cdot \|_{n, {\infty} } )} du +
%C_0 MK $$ 
%$$=2C_0  \sqrt {\log (2p)\log (2n)} MK_X \log n + C_0 MK_X . $$
%\hfill $\sqcup \mkern -12mu \sqcap$ 

As a consequence, we obtain a result which is in \cite{rudelson2011reconstruction}. 
 It suffices to combine Theorem \ref{square.theorem} with Lemma \ref{Rudelson.lemma}.
\begin{theorem} \label{l1.square.theorem}
For all $t >0$
$$\PP \biggl (\sup_{\| \beta \|_1 \le M,\ \| f_{\beta} \| \le 1  } \biggl | \| f_{\beta} \|_n^2 - \| f_{\beta}  \|^2 \biggr | /c_1 \ge
 MK_X \sqrt { \log p \log^{3} n  +t \over  n  }   $$ $$ + M^2 K_X^2 {    \log p \log^3 n  +t  \over n}
  \biggr )  \le \exp[-t]  .$$
%  as well as
%  $$\PP \biggl (\sup_{\| \beta \|_1 \le M, \| f_{\beta} \| \le 1  }\ \sup_{\| \tilde \beta \|_1 \le M, \| f_{\tilde \beta} \| \le 1  }
%   \biggl | (P_n - P) f_{\beta} f_{\tilde \beta}  \biggr | /c_1 \ge
% MK_X \sqrt { \log p \log^{3} n  +t  \over  n  }   $$ $$ + M^2 K_X^2{      \log p \log^3 n  +t   \over n}
%  \biggr )  \le \exp[-t]  .$$

\end{theorem}

Theorem \ref{l1.square.theorem} has very useful applications, in particular to $\ell_1$-regularization
or to exact recovery using basis pursuit (\cite{chen1998atomic})
where results often rely on bounds for compatibility constants (\cite{vandeG07}, \cite{vdG:2009}) or restricted eigenvalues
(\cite{bickel2009sal}). This is elaborated
upon in \cite{rudelson2011reconstruction}. 

Theorem \ref{l1.square.theorem} can be applied also to obtain a uniform bound over all
subspaces. 
Define the minimal eigenvalue $\Lambda_{\rm min}^2  := \min_{\| \beta \|_2 \le 1 } \| f_{\beta} \|^2 $.

\begin{theorem} \label{l0.square.theorem} Suppose $\Lambda_{\rm min} > 0 $. 
Define for $S \subset \{ 1 , \ldots , p \}$, $\beta_{j,S} = \beta_j {\rm l} \{ j \in S\} $,
$j=1 , \ldots , p $.
For all $t>0$
\begin{equation}\label{l0.equation}
\PP \biggl (\exists \ s:\ \sup_{|S|=s} \sup_{\| f_{\beta_S}  \| \le 1  } \biggl | \| f_{\beta_S} \|_n^2 - \| f_{\beta_S}  \|^2 \biggr |  /c_1\ge
{  K_X \over \Lambda_{\rm min}}  \sqrt {s  \log p \log^3 n +st   \over  n  }    
\end{equation}
 $$+
 {  K_X^2 \over \Lambda_{\rm min}^2 }  {  s  \log p \log^3 n  +st \over n}
  \biggr )  \le \exp[-t]  . $$
% as well as
%$$\PP \biggl (\sup_{\| f_{\beta_S}  \| \le 1  } \sup_{\| f_{\tilde \beta_S} \| \le 1 } \biggl | (P_n - P) f_{\beta_S} f_{\tilde \beta_S} \biggr |  /c_1\ge
%{  K_X \over \Lambda_{\rm min}}  \sqrt {s  \log p \log^3 n  +st   \over  n  }    $$ $$+
% {  K_X^2 \over \Lambda_{\rm min}^2 }  {  s  \log p \log^3 n  +st \over n}
%  \biggr )  \le \exp[-t]   .$$
\end{theorem}

The next theorem is a direct application of a Bernstein type inequality for random matrices as given in 
\cite{ahlswede2002strong} (see also Theorem 3 in \cite{koltchinskii2011remark}). It shows that in  Theorem \ref{l0.square.theorem} the $\log^3 n$-term can be omitted when one considers a fixed set $S$
instead of requiring a result uniform in $S$.

\begin{theorem}\label{l3.square.theorem}Suppose $\Lambda_{\rm min} > 0 $. 
For all $t>0$
\begin{equation}\label{l3.equation}
\PP \biggl ( \sup_{\| f_{\beta}  \| \le 1  } \biggl | \| f_{\beta} \|_n^2 - \| f_{\beta}  \|^2 \biggr |  /c_1\ge
{  K_X \over \Lambda_{\rm min}}  \sqrt {p  \log p  +pt   \over  n  }    
\end{equation}
 $$+
 {  K_X^2 \over \Lambda_{\rm min}^2 }  {  p  \log p   +pt \over n}
  \biggr )  \le \exp[-t]  . $$
% as well as
%$$\PP \biggl (\sup_{\| f_{\beta}  \| \le 1  } \sup_{\| f_{\tilde \beta} \| \le 1 } \biggl | (P_n - P) f_{\beta} f_{\tilde \beta} \biggr |  /c_1\ge
%{  K_X \over \Lambda_{\rm min}}  \sqrt {p \log p   +pt   \over  n  }    $$ $$+
% {  K_X^2 \over \Lambda_{\rm min}^2 }  {  p \log p \log n  +pt \over n}
%  \biggr )  \le \exp[-t]   .$$

\end{theorem} 

\begin{remark} Let us briefly indicate how this compares to an isotropic case.
Following an idea of  \cite{loh2012} (see also Lemma 1 in \cite{Nickl:2013}) one can show that the supremum over all $\| f_{\beta} \| \le 1$ can in fact 
be replaced by a maximum over a finite class: 
$$\sup_{ \| f_{\beta} \| \le 1 } \biggl | \| f_{\beta} \|_n^2 - \| f_{\beta} \|^2 \biggr | \le c_1
\max_{j \in \{ 1 , \ldots , N\} } \biggl | (P_n - P ) f_{\beta_j }^2 \biggr | , $$
where $\| f_{\beta_j} \| \le 1$ for all $j=1, \ldots , N$ and where $\log N \le c_0^2 p $.
We can now proceed by invoking the union
bound for the maximum. An isotropy assumption then leads to good results.
 We assume sub-Gaussianity
of the vectors $\{ X_i \}$, meaning that each $f_{\beta} (X_i) $ is sub-Gaussian: there is a constant $K_1$ such that
for all $\| f_{\beta} \| \le 1$ and all $i$ it holds that
$\| f_{\beta} (X_i) \|_{\Psi_2 } \le K_1$. Then by Bernstein's inequality, for all all $\| f_{\beta} \| \le 1 $ and all $t >0$ 
$$ \PP \biggr (\biggl | \| f_{\beta} \|_n^2 - \| f_{\beta} \|^2 \biggr |  / C_1  \ge K_1^2 \sqrt {t / n} + K_1^2 t/n \biggl ) \le \exp [-t] . $$
The union bound together with the above reduction then gives for all $t >0$
$$\PP \biggl ( \sup_{ \| f_{\beta} \| \le 1 } \biggl | \| f_{\beta} \|_n^2 - \| f_{\beta} \|^2 \biggr | /(c_1 C_1) \ge
K_1^2  \sqrt { c_0^2 p + t \over n} $$ $$+ K_1^2 {c_0^2 p +t  \over n }\biggr ) 
\le \exp[-t] . $$
The latter result is a ``true" deviation inequality: the deviation from the bound $\asymp \sqrt {p/n} $ for the mean does not involve this bound, i.e., there is no $p$ in front of $t$ inside the probability. This in
contrast to the result (\ref{l0.equation}) in Theorem \ref{l3.square.theorem}.

\end{remark}

%\begin{remark} Note that the $(\log p \log^3 n )$-term of Theorem \ref{l1.square.theorem} is replaced
%by a $(\log^3 p \log n )$-term. This (admittedly very small) improvement follows from the fact that for small values
%of the integrant in the entropy integral one can replace the entropy estimate of Lemma \ref{Rudelson.lemma}
%by another - more standard - one for $p$-dimensional spaces,  as is done in \cite{rudelson2011reconstruction}.
%It is to us not clear whether the $\log p$- and/or $\log n$-terms can be improved in Theorem \ref{l0.square.theorem}.

%\end{remark}

\begin{remark} One may wonder why the minimal eigenvalue is playing a role
in the result of Theorems \ref{l0.square.theorem} and  \ref{l3.square.theorem}. Of course, as far as conditions on $L_2 (P)$-norms are concerned
one may orthogonalize the variables. However, after orthogonalization 
the sup-norm of the variables can be quite large. 
\end{remark}

The following lemma improves Theorem \ref{subGaussianproduct.theorem} in the linear case.

\begin{lemma} \label{Gaus.lemma} Suppose that ${\bf Y} := \{ Y_1 , \ldots , Y_n \}$ is uniformly sub-Gaussian
with constant $K_0$ (see Definition \ref{conditional-subGaussian.definition}). 
Then for all $t >0$
$$ \PP \biggl ( \sup_{\| f_{\beta} \|_1 \le 1 } | (P_n - P){\bf Y}  f_{\beta} |/ c_2  \ge{K_0 K_X \over  \Lambda_{\rm min}}\sqrt { pt \over n }  \biggr ) \le 2 \exp [-t] . $$

\end{lemma}

%We now apply Theorem \ref{subGaussianproduct.theorem} (products with a
%sub-Gaussian random variable) to our special case of linear functions.
%
%\begin{lemma}\label{Gaus2.lemma}
%Suppose ${\bf Y}$ is conditionally sub-Gaussian given ${\bf X}$ with constant $K_0$.
% Consider values of $t\ge 4$ and $n$ such that
%$$ \sqrt {K_X^2 p \log^3 p \log n (1+t)  \over   n  \Lambda_{\rm min}^2 }  +  {K_X^2  p \log^3 p  \log n (1+t) \over   n
%\Lambda_{ \rm min}^2}  
%  \le {1 \over c_0} $$
%For these values
%$$\PP \biggl (\sup_{\| f_{\beta} \| \le 1 }  | (P_n - P) Y f_{\beta}   |/ c_1 \ge
%K_0  \sqrt {p \log p   \over   n \Lambda_{\rm min}^2  }  +  K_0  \sqrt {t \over n} 
%  \biggr )  \le 8 \exp[-t]. $$
%\end{lemma}

To avoid too involved expressions, we from now on will use order symbols.
Then, the results needed for the next section can be summarized as follows. 

\begin{summary}\label{linear.summary}
Suppose that ${\bf Y} := \{ Y_1 , \ldots , Y_n \}$ is uniformly sub-Gaussian
with constant $K_0$, that 
$ \max_{i,j} | X_{i,j} | \le K_X $, $\Lambda_{\rm min} > 0 $ and that $\delta_n = o(1)$, where
$$ \delta_n^2 := { K_X^2(1 + K_0^2)  p \log p  \over n \Lambda_{\rm min}^2 } .$$
Then uniformly in $\| f_{\beta} \| \le 1 $, $\| f_{\tilde \beta} \| \le 1 $ is holds that
$$\biggr |  \| f_{\beta} \|_n^2 - \| f_{\beta } \|^2\biggr  | ={\mathcal O}_{\PP} (\delta_n) , \ 
\biggl | (P_n - P) ({\bf Y} - f_{\tilde \beta} ) f_{\beta}\biggr  |  ={\mathcal O}_{\PP}  (\delta_n )  . $$
 \end{summary}

\section{Least squares when the model is wrong}\label{least.squares.section}
In this section we examine a $p$-dimensional linear model with $p$ moderately large, and the
least squares estimator.
The observations are $\{ (X_i , Y_i) \}_{i=1}^n $, independent, and with
$X_i \in {\cal X}$ and $Y_i \in \R $ ($i=1 , \ldots , n $). 
Let $\{ \psi_j \}_{j=1}^p $ be a given dictionary of functions on ${\cal X}$.
We write $f_{\beta} (\cdot) := \sum_{j=1}^p \beta_j \psi_j (\cdot) $, $\beta \in \R^p$.

The least squares estimator is
$$ \hat f := \arg \min_{f_{\beta} } \sum_{i=1}^n ( Y_i - f_{\beta} (X_i) )^2 .$$
Let $f^0 (X_i) := \EE( Y_i \vert X_i)$ be the conditional
expectation of $Y_i$ given $X_i$,  $i=1 , \ldots , n $.
The projection in $L_2 (P)$ of $f^0$ on the linear space $\{ f_{\beta}: \beta \in \R^p \}$ is written as $f^*$.
We want to show convergence of $\hat f$ to $f^*$. Because we know little about the higher order
moments of $f^*$ (only the second moment is under control as $\| f^* \| \le \| {\bf Y} \|$)
the situation is a little more delicate than in the usual regression context (where
$\| f^* - f^0 \| $ is small).  This is where uniform convergence
of $\| \cdot \|_n$ to $\| \cdot \|$ comes in. 

\begin{lemma} \label{basic.lemma} Let $0 < \delta_n < 1/2 $. On the set
$${\cal T} := \bigg \{ \sup_{\| f_{\beta} \| \le 1 }  \biggl |  \| f_{\beta }\|_n^2 - \| f_{\beta} \|^2 \biggr | \le
\delta_n, \   \sup_{\| f_{\beta} \| \le 1, \ \| f_{\tilde \beta} \| \le 1  }  \biggl | 2 (P_n -P) ({\bf Y}- f_{\tilde \beta})  f_{\beta}\biggr  | \le \delta_n  \biggr \} ,$$
it holds that
$$\| \hat f - f^* \| \le 2 \delta_n . $$

\end{lemma}

To handle the set ${\cal T}$ given in the above lemma, we invoke Summary \ref{linear.summary}.
To this end, define the matrix $ \Sigma := P \psi^T \psi $
and let $\Lambda_{\rm min}^2 $ be the smallest eigenvalue of $\Sigma$. 

\begin{theorem}\label{theoretical.norm.theorem}
Suppose that ${\bf Y} := \{ Y_1 , \ldots , Y_n \}$ is uniformly sub-Gaussian
with constant $K_0$ (see Definition \ref{conditional-subGaussian.definition}), 
that $ \max_{i,j} | X_{i,j} | \le K_X$, 
$\Lambda_{\rm min} >0 $, and that $\delta_n = o(1)$ where
$$ \delta_n^2 := { K_X^2(1 + K_0^2)  p \log p  \over n \Lambda_{\rm min}^2 } .$$
 Then
$$ \| \hat f - f^* \| = {\mathcal O}_{\PP}  ( \delta_n   ) . $$
Moreover
$$ \biggl | \| {\bf Y}- \hat f \|_n^2 - \| {\bf Y} -f^* \|^2 \biggr | = 
{\mathcal O}_{\PP}  (\delta_n     ) . $$

\end{theorem}

      In view of the uniformity in Summary \ref{linear.summary} we can formulate an extension. Such an
      extension will 
      be useful in the next section. 
      Recall the notation: for a set $S \subset \{ 1 , \ldots , p \}$ and $\beta \in \R^p$
      $$\beta_{j,S} := \beta_j \{ j \in S \} , \ j=1 , \ldots , p .  $$
      Consider, for any set $S \subset \{ 1 , \ldots , p \}$, the projection $f_S^*$ of $f^0$ on the $|S|$-dimensional space
    ${\cal F}_S := \{ f_{\beta_S } : \ \beta \in \R^p \} $ and the corresponding least squares estimator 
    $$\hat f_S: = \arg \min_{f_{\beta_S} } \| {\bf Y} - f_{\beta_S} \|_n . $$

    \begin{theorem}\label{theoretical.norm2.theorem} Assume the conditions of Theorem \ref{theoretical.norm.theorem} and let $\delta_n$ be defined as there.
    Then uniformly in all $S$,
    $$ \| \hat f_S - f_S^* \| = {\mathcal O}_{\PP}  ( \delta_n   ) . $$
Moreover
$$ \biggl | \| {\bf Y}- \hat f_S \|_n^2 - \| {\bf Y}-f_S^* \|^2 \biggr | = 
{\mathcal O}_{\PP}  (\delta_n  ) . $$
 \end{theorem}

   \section{Application to DAG's}\label{DAG.section}
   
   Let $X$ be a $n \times p $ matrix with i.i.d. rows. We throughout this section assume $p \le n$.
   The $i$-th row is denoted by $X_i:= (X_{i,1} , \ldots , X_{i,p} )$ ($i=1 , \ldots , n$). 
      The distribution of a row, say  $X_1$, is denoted by $P$. 
   
   We assume a directed acyclic graph (DAG) structure. Namely,  we assume the structural equations model defined
   as follows.
   
   \begin{definition} \label{sem.definition} We say that $X_1 \in \R^p$ satisfies the non-linear Gaussian structural
   equations model if for some permutation $\pi^0$ of $\{ 1 , \ldots , p \}$ and
   for some functions $f_j^0 : \ \R^{j-1} \rightarrow \R$
  $$ X_{1,\pi_j^0} = f_j^0 (   X_{1, \pi_1^0 } , \ldots X_{1, \pi_{j-1}^0} ) +
   \varepsilon_{1,\pi_j^0}  ,\ j=1 , \ldots , p , $$
   where $\{ \varepsilon_{1, 1 }, \ldots , \varepsilon_{1,p} \}$ are independent and where
   for $j=1 , \ldots , p $ the random variable
   $\varepsilon_{1,\pi_j^0}\sim  {\cal N} (0, \sigma_{\pi_j^0}^2 )$ is independent of
   $( X_{1, \pi_1^0 } , \ldots ,  X_{1, \pi_{j-1}^0})$. The latter set is to be
   understood as the empty set when $j=1$. 
   
  \end{definition}
  
    We let $\Pi_0$ be the set of permutations $\pi_0$ for which Definition \ref{sem.definition} holds.
    The case of interest is the one where
    $f_j^0 (  X_{1, \pi_1^0 } , \ldots ,  X_{1, \pi_{j-1}})= f_j^* ( \{ X_{1,k} \}_{k \not= j } \} ) $,
    $j=1 , \ldots , p$, with 
    with $\{ f_j^*\}$ and hence
    and $\{ \sigma_j^2\} $ not depending on $\pi_0$.
   Our aim is to find a member from  $\Pi^0$ based on the data
   $X_1 , \ldots , X_n$.

   \subsection{Some notation}\label{notation.section}
   
   We consider a given class ${\cal F}_0$ of functions $f_0 : \ \R \rightarrow \R$.

    Let ${\cal F}_1 := \emptyset$ and for $j=2 , \ldots , p $
    $${\cal F}_j  := \biggl \{ f(x_1 , \ldots , x_{j-1} ) = \sum_{k=1}^{j-1} f_{k,j} (x_k) :  $$ $$ \ \ \ \ \ \ \ 
    (x_1 , \ldots , x_{j-1} ) \in \R^{j-1} , \ f_{k,j} \in {\cal F}_0 \ \forall \ j,k \biggr \} . $$
    Let $\Pi$ be the set of all permutations of $\{ 1 , \ldots , p \}$. 
   Write for each permutation $\pi \in \Pi $ and each $j$, for $f_j \in {\cal F}_j$,
   $$ f_j (X_i , \pi)  :=  f_j ( X_{i, \pi_1} , \ldots , X_{i, \pi_{j-1}}  ), \ i=1 , \ldots , n . $$
   Define
$$f_j^* (\pi) :=  \arg \min  \biggl \{  \| {\bf X}_{ \pi_j} -   f_j( \pi)  \|^2 :  \ 
   f_j \in {\cal F}_j \} .  $$
where for $j=1$, we take $f_1^* (\pi ) = 0 $,
and where
$$ \| {\bf X} _{\pi_j} -  f_j(\pi) \|^2:= 
 \EE\biggl  ( X_{1,  \pi_j} -
 f_j ( X_1, \pi )\biggr  )^2 . $$

We further define
   $$\sigma_j^2 (\pi) := \| {\bf X}_{ \pi_j} - f_j^* (\pi) \|^2 , \ j=1 , \ldots , p .  $$
   
%   We finally let
%   $${\cal F}:= \biggl \{ f = \sum_{k=1}^p f_k (x_k) , \ f_k \in {\cal F}_0 \ \forall k \biggr \} . $$
%   
    
   \subsection{Identifiability}
   
   In order to be able to estimate a correct permutation one needs to assume that the
wrong permutations can be detected. 

   \begin{condition}\label{identifiable.condition} (Identifiability condition). For some constant $\xi >0$, 
   $$ \inf_{\pi \notin \Pi_0, \ \pi_0 \in \Pi_0 } {1 \over p}  \sum_{j=1}^p \log \biggl ( { \sigma_j ( \pi)  \over  \sigma_j (\pi^0) } \biggr ) >   \xi  .  $$
   \end{condition}
   
   This condition is discussed in \cite{PetersB}. The linear Gaussian structural equations model has
   $\Pi_0= \Pi$, i.e. any permutation is correct. In the non-linear case, we think
   of the situation where, unlike the linear case, the parental dependence is the same
   for all $\pi_0 \in \Pi_0$, say
   $f_j^0 (X_{1, \pi_1^0 } , \ldots , X_{1 , \pi_{j-1}^0 } ) :=
   f_j^* (\{ X_{1,k}  \}_{k \not= j } )$ ($j=1 , \ldots , p$),
   and hence also the residual variances $\sigma_j^2$, $j=1 , \ldots , p $
   do not depend on $\pi_0$.   
   The identifiably condition then requires that choosing $\pi \notin \Pi_0$ will
   give on average too large residual variances.  If the model is misspecified, Condition
   \ref{identifiable.condition} is to be seen as assuming robustness to the bias that misspecification
   introduces. In an asymptotic formulation, it suffices to assume identifiability at the truth:
  $ \inf_{\pi \notin \Pi_0 }  \sum_{j=1}^p \log  ( { \sigma_j (\pi) / \sigma_j }  ) /p > \xi_0$
   with $1/ \xi_0= {\mathcal O} (1)$ together with a vanishing bias: 
   $\sup_{\pi_0 \in \Pi_0} \sum_{j=1}^p \log (\sigma_j (\pi_0) / \sigma_j )/ p \rightarrow 0 $.
  One may consider choosing a model with low complexity (large bias)
   because $\pi_0$ is the parameter of interest here. The estimation of $f_j^0$
   ($j=1 , \ldots , p$) can then follow in a second step using a standard (nonparametric)
   regression estimator and the estimated permutation.    
   
   \subsection{The estimator}
To describe the estimator of $\pi^0$ we introduce empirical counterparts of the
quantities given above. For each $j$ and $\pi$ we write
$$\| {\bf X}_{\pi_j} - f_j (\pi) \|_n^2 :=  {1 \over n} \sum_{i=1}^n  \biggl  ( X_{i,  \pi_j} -
f_j( X_i, \pi   )\biggr  )^2 . $$
We let $\hat f_j (\pi)$ be the least squares estimator
$$\hat f_j (\pi) := \arg \min \biggl \{  \| {\bf X} _{\pi_j} - f_j (\pi) \|_n :  \ 
  f_j \in {\cal F}_j  \biggr \}   $$
and take the normalized residual sum of squares
$$\hat \sigma_j (\pi) := \| {\bf X}_{\pi_j}- \hat f_j (\pi) \|_n^2 $$
as estimator of $\sigma_j^2 (\pi)$. 
We then let
\begin{equation}\label{hatpi.equation}
\hat \pi :\in  \arg \min_{\pi \in \Pi } \sum_{j=1}^p \log  \hat \sigma_j^2 (\pi)  . 
\end{equation}

\subsection{Consistency}

Let $H( \cdot , {\cal F}_0 , \| \cdot \|_{\infty}) $ be the entropy of ${\cal F}_0$ endowed with
supremum norm.

\begin{theorem} \label{DAG0.theorem} Suppose the non-linear Gaussian structural equations model
(see Definition \ref{sem.definition}) with $ \max_{1 \le j \le p } \sigma_j^2 = {\mathcal O} (1)$
and
$ \max_{1 \le j \le p } \| f_j^0 \|_{\infty} = {\mathcal O} (1) $. 
Assume Condition \ref{identifiable.condition} (the identifiability condition)
with $1/ \xi = {\mathcal O}  (1)$.
Assume moreover that ${\cal F}_0$ is a convex class and that 
one of the following 5 cases hold of Theorem \ref{Sobolev.theorem} for the collection ${\cal F}:=
\{ f = \sum_{k=1}^p f_k (x_k) , \ f_k \in {\cal F}_0 \ \forall k  \} $:\\
{\bf Case 1.} Case 1  holds and $p^3 /n = o(1)$,\\
{\bf Case 2.} Case 2 holds and $p^{3-(1- \alpha)^2} /n = o(1)$,\\
{\bf Case 3.} Case 3 holds and $p^{2}\log^4 n  /n = o(1)$,\\
{\bf Case 4.} Case 4 holds and $pN^2 \log n  /n = o(1)$,\\
{\bf Case 5.} Case 5 holds, $p^{1+ 4 \alpha} \log n  /n = o(1)$.  \\
Then $
\PP ( \hat \pi \notin \Pi_0 ) \rightarrow 0 $. 
\end{theorem}

We recall Remark \ref{local.remark}: the conditions on ${\cal F}$ may be weakened to local
versions.

\section{Conclusion} \label{conclusion.section}
In this paper we summarized  some results for the uniform
convergence of empirical norms and the extension
to empirical inner products. 

For statistical theory the results are very useful.
In \cite{bartlett2012} one can find an application to $\ell_1$-restricted
regression for the case of random design and \cite{rudelson2011reconstruction} focuses
on the restricted isometry property and restricted eigenvalues. 
We have given the application to order estimation in directed acyclic graphs (DAG's). 
We omitted important computational issues and further discussions for this special
case as it is beyond the scope of the paper. For more details we refer to \cite{PetersB}.

The results can also
be applied to generalize the results in \cite{vdGP:2013} for DAG's to the linear non-Gaussian case, in particular to
anisotropic distributions. A generalization to 
to isotropic distributions (e.g. sub-Gaussian distributions) is possible but perhaps less relevant as in many statistical applications isotropy is not very natural or stable (for DAG's sub-Gaussianity can hold when the linear model is exactly true but
it is not clear what happens when the model is only approximately
linear). 

A further application is the estimation of a precision matrix for non-Gaussian data.
We mention that such an approach is used in \cite{vdGBR:2013} to construct confidence
intervals for a single parameter.
Here, a Lasso is used for estimating a Fisher-information matrix. The estimator is based on 
empirical projections
and also the function to be estimated is a theoretical projection
as in Section \ref{least.squares.section}. 
In the context of confidence intervals in $\ell_2$,  the uniform
convergence may generalize the (sub-)Gaussian case considered in \cite{Nickl:2013}.
Another application of uniform convergence, this time 
for additive models (\cite{MullervdGeer13},
\cite{vdGMammen13}), 
was
briefly indicated in Subsection \ref{smooth.section}.

% Let 
%$$ \alpha := {1 \over b_* + \| \hat f - f^0 \| } ,$$
%and $\hat f_{\alpha} :=  \alpha \hat f + (1-\alpha ) f^*$.
%Then we know that
%$$\| \hat f_{\alpha} - f^0 \|_n \le 2 P_n \varepsilon (\hat f_{\alpha } - f^* ) +
%\| f^* - f_0 \|_n . $$
%We also know that
%$$\| \hat f_{\alpha} - f^* \| = \alpha \| \hat f - f^* \| \le b_* . $$
%So
%$$\| \hat f_{\alpha} - f^0 \|^2 \le 2 b_*^2 . $$
%So 
%$$\| \hat f_{\alpha} - f^0 \|^2  \le \| \hat f_{\alpha} - f^0 \|_n^2 + 2 b_*^2 \delta_n \le
%\tilde \delta_n \| \hat f_{\alpha} - f^* \| + b_*^2  + 3 b_*^2 \delta_n  $$
%$$ \le  \tilde \delta_n^2 + \| \hat f_{\alpha} - f^* \|^2/2 +  b_*^2  + 3 b_*^2 \delta_n .$$
%So
%$$\| \hat f_{\alpha} - f^* \|^2 \le 2 \tilde \delta_n^2 + 6 b_*^2 \delta_n . $$
%It follows that
%$$\| \hat f_{\alpha} - f^* \| \le \sqrt 2 \tilde \delta_n + \sqrt 6 b_* \sqrt {\delta_n} . $$
%But then
%$$\| \hat f - f^* \| \le \sqrt 2 \tilde \delta_n  b_* + \sqrt 6 b_*^2 \sqrt {\delta_n} +
%\sqrt 2 \tilde \delta_n  \| \hat f - f_* \| + \sqrt 6 b_* \sqrt {\delta_n} \| \hat f - f^* \| . $$
%So then we know that $\| \hat f - f^* \| \le 1 $.
%But we also know that
%$$ \| \hat f - f^* \|_n^2 \le \| f^* - f^0 \|_n^2 + \bar \delta_n^2 . $$ 
%So now,
%$$ \| \hat f - f^* \|^2 \le \| f^* - f_0 \|^2 + 2 \delta_n + \bar \delta_n^2 . $$

\section{Technical tools}\label{technical.section}

\subsection{Symmetrization}

Define
$${\bf Z} ({\cal F} ) := \sup_{f \in {\cal F} } \biggl | (P_n - P) f \biggr | . $$
Let moreover $\epsilon_1 , \ldots , \epsilon_n$  be a Rademacher sequence (that is,
$\epsilon_1 , \ldots , \epsilon_n $ are independent random variables taking the values
$+1$ or $-1$ each with probability ${1 \over 2} $) independent of
$X_1 , \ldots , X_n$,  and define
$${\bf Z}^{\epsilon} ({\cal F} ) := \sum_{i=1}^n \epsilon_i f(X_i) / n . $$

\begin{theorem}\label{symmetrization.theorem} (see e.g.\ \cite{vanderVaart:96}). It holds that
$$ \EE {\bf Z} ({\cal F} ) \le 2 \EE {\bf Z}^{\epsilon} ({\cal F} ) . $$
\end{theorem}

\begin{theorem}\label{symmetrization2.theorem} (see \cite{pollard1984convergence}).
Let $R:= \sup_{f \in {\cal F} }\| f \| $. For $t \ge 4  $,
$$ \PP ( {\bf Z} ({\cal F}) \ge 4R \sqrt {2t /n}) \le 4 \PP ({\bf Z}^{\varepsilon} ({\cal F} ) \ge R \sqrt {2t/n} ) . $$
\end{theorem} 

\subsection{Dudley's theorem}
%For $k=1, 2$, we let $\| \cdot  \|_{\Psi_k}$ be the Orlicz norm corresponding to the function $\Psi_k (z) := \exp [|z|^k] -1$,
%$z \in \R$. If for a random variable $Z \in \R$, $\| Z \|_{\Psi_2} < \infty$ we call $Z$ sub-Gaussian.
%If $\| Z \|_{\Psi_1} < \infty$ we call it sub-exponential.  

Dudley's theorem is originally for Gaussian processes (see \cite{dudley1967sizes}). The extension to
sub-Gaussian random variables and Rademacher averages is rather straightforward. We summarize these
in our context 
in Theorem \ref{Dudley.theorem} below.

%Let ${\cal H}  (\cdot , {\cal F} , \| \cdot \| ) $ denote the entropy
%of ${\cal F}$ equipped with the metric induced by the $L_2 (P)$-norm $\| \cdot \|$ and 
%let $R := \sup_{f \in {\cal F}}  \|f \| $ be the radius of ${\cal F}$.

Let ${\cal H}  (\cdot , {\cal F} , \| \cdot \|_n ) $ denote the entropy
of ${\cal F}$ equipped with the metric induced by the empirical norm $\| \cdot \|_n$.
and let $\hat R$ be the random radius $ \hat R := \sup_{f \in {\cal F} } \| f \|_n $.

%\begin{theorem}\label{Dudley0.theorem} (sub-Gaussian case).
%Suppose that for all $f , \tilde f  \in {\cal F}$
%$$\| (P_n - P) (f - \tilde f)  \|_{\Psi_2} \le \|\tilde f  - f  \|/ \sqrt n  . $$
%Then
%$$ \EE {\bf Z} ({\cal F} ) \le C_0 R  \int_0^1 \sqrt {{\cal H} ( u R , {\cal F} , \| \cdot \| ) } du / \sqrt n .$$
%\end{theorem}

\begin{theorem}\label{Dudley.theorem} (Rademacher averages). We have
$$ \EE {\bf Z}^{\epsilon} ({\cal F} ) \le C_0 \inf_{\delta > 0 } \EE \biggl [ \hat R  \int_{\delta}^1 \sqrt {{\cal H} ( u \hat R , {\cal F} , \| \cdot \|_n ) } du / \sqrt n  + \delta \hat R\biggr ] .$$
\end{theorem}

\subsection{Deviation inequalities} 

We present two deviation inequalities, for the bounded case and the sub-Gaussian case.

\begin{theorem} \label{bounded.theorem} (see \cite{Talagrand:95}, \cite{Massart:00a}). Suppose that for some constants $R$ and $K$.
$$ \sup_{f \in {\cal F} } \| f \| \le R , \ \sup_{f \in {\cal F}} \| f \|_{\infty} \le K . $$
Then for all $t>0$
$$ \PP \biggl ( {\bf Z} ({\cal F}) / C_1 \ge   \EE {\bf Z} ({\cal F}) + R \sqrt {t / n} + K t / n  \biggr ) \le \exp[-t] . $$
\end{theorem}

\begin{theorem}\label{subGaussian.theorem} (see \cite{Massart:00a}).
Let ${\bf X} := (X_1 , \ldots , X_n)$. Conditionally on ${\bf X}$, for
all $t > 0$, 
$$ \PP \biggl ( {\bf Z}^{\epsilon} ({\cal F}) / C_1 \ge  \EE(  {\bf Z}^{\epsilon} ({\cal F} ) \vert {\bf X} )
 + \hat R  \sqrt {t / n} \biggl \vert {\bf X}  \biggr ) \le \exp[-t] , $$
where $\hat R := \sup_{f \in {\cal F}} \| f \|_n $. 
\end{theorem}

%\begin{theorem}\label{subGaussian.theorem} (see \cite{viens2007supremum}).
%Suppose that for all $f , \tilde f \in {\cal F} $
%$$\| (P_n - P) (f - \tilde f)  \|_{\Psi_2} \le \|\tilde f  - f  \| / \sqrt n . $$
%Then for all $t > 0$, 
%$$ \PP \biggl ( {\bf Z} ({\cal F}) / C_1 \ge  C_0 R  \int_0^1 \sqrt {{\cal H} ( u R , {\cal F} , \| \cdot \| ) } du / \sqrt n
% + R  \sqrt {t / n}  \biggr ) \le \exp[-t] , $$
%where $R := \sup_{f \in {\cal F}} \| f \| $. 
%\end{theorem}

%Let $F (\cdot) := \sup_{f \in {\cal F}} | f (\cdot ) | $ be the envelope function of $( {\cal F}, \| \cdot \| $.
%
%\begin{theorem} \label{sub-exponential.theorem} (see \cite{Adamczak:08}). Suppose that for some constants $R$ and $K$
%$$ \sup_{f \in {\cal F} } \| f \| \le R , \  \max_{i} \biggl \| | F(X_i) | \biggr  \|_{\Psi_1} \le K   . $$
%Then for all $t > 0$, 
%$$ \PP \biggl ( {\bf Z} ({\cal F}) /C_1 \ge   \EE {\bf Z} ({\cal F}) +  R \sqrt {t / n} +  K t \log n  / n ) \biggr ) \le \exp[-t] . $$ 
%\end{theorem}
%
%The factor $\log n$ can be removed in certain cases, see \cite{vandeGeerLederer:13}. 
%

\section{Proofs} \label{proofs.section}

\subsection{Proofs for Section \ref{empirical.norm.section}}

 Theorem \ref{square.theorem} follows from \cite{guedon2007subspaces}.
We present a proof for completeness and to facilitate the extension to products of functions.

{\bf Proof of Theorem \ref{square.theorem}.} We consider the symmetrized process
$$P_n^{\epsilon} f^2 := \sum_{i=1}^n \epsilon_i f^2 (X_i) /n, $$
with $\epsilon_1 , \ldots , \epsilon_n$ a Rademacher sequence independent of $X_1 , \ldots , X_n $,
and then apply Dudley's theorem, see Theorem \ref{Dudley.theorem}.
Note that for two functions $f $ and $\tilde f$ in the class ${\cal F}$
$$ \| f^2 - \tilde f^2  \|_n  \le
\| f+ \tilde f \|_n \| f - \tilde f \|_{n, \infty} \le2  \hat R \| f - \tilde f \|_{n, \infty} .$$
It follows that
$${\cal H} ( u , {\cal F}^2 , \| \cdot \|_n ) \le {\cal H} ( u/ (2 \hat R)  , {\cal F} , \| \cdot \|_{n, {\infty} } ) , \ u > 0 . $$
Hence
$$\int_{\delta}^{1} \sqrt {{\cal H} (u \hat R K , {\cal F}^2 , \| \cdot \|_n }du
 \le \int_{\delta}^{1} \sqrt {{\cal H} ( uK /2   , {\cal F} , \| \cdot \|_{n, {\infty} } )} du .$$
Here we used that $\|  f^2 \|_n \le  \hat R K$. 
So by Theorem \ref{Dudley.theorem}
$$\EE \biggl ( \sup_{f \in {\cal F} } P_n^{\epsilon} f^2 \biggr ) \le 
C_0 \inf_{\delta >0}  \EE \biggl [ \hat R K \int_{\delta}^{1} \sqrt { {\cal H} (uK/2 , {\cal F} , \| \cdot \|_{n, {\infty} } )} du
/ \sqrt {n} 
+ \delta \hat R K   \biggr ]  $$
$$ \le  J_{\infty} (K, {\cal F} ) \sqrt { \EE \hat R^2 }/\sqrt n   . $$
But then by Theorem \ref{symmetrization.theorem}
\begin {equation}\label{first-bound.equation}
 \EE \biggl (\sup_{f \in {\cal F} } \biggl |  \| f \|_n^2 - \| f \|^2 \biggr | \biggr ) \le
2 J_{\infty} (K, {\cal F} ) \sqrt { \EE \hat R^2 } / \sqrt {n}. 
\end{equation}
This leads to the by-product of the theorem: the inequality
$$ \EE \hat R^2 \le R^2 + 2 J_{\infty} (K, {\cal F} ) \sqrt { \EE \hat R^2 }/ \sqrt n $$
gives
$$\sqrt  {\EE \hat R^2 } \le R + 2 J_{\infty} (K, {\cal F} )/\sqrt n  . $$
Insert this in (\ref{first-bound.equation}) to find
$$ \EE \biggl (\sup_{f \in {\cal F} } \biggl |  \| f \|_n^2 - \| f \|^2 \biggr | \biggr ) \le 2 J_{\infty}  (K, {\cal F}) R /\sqrt n +
4 J_{\infty}^2 (K,  {\cal F })/n  . $$ 
We now apply Theorem \ref{bounded.theorem}. We have
$$\sup_{f \in {\cal F}} \| f^2 \| \le RK , \ \sup_{f \in {\cal F} } \| f^2 \|_{\infty} \le K^2 . $$
Hence, inserting the just obtained bound for the expectation, for all $t > 0$ 
$$\PP \biggl (\sup_{f \in {\cal F}} \biggl | \| f \|_n^2 - \| f \|^2 \biggr | / C_1 \ge
{2  R J_{\infty}  (K , {\cal F} ) + RK \sqrt t  \over  \sqrt n  } + { 4  J_{\infty}^2 (K, {\cal F}) + K^2 t \over n}    \biggr ) $$ $$ \le \exp[-t] . $$

\hfill $\sqcup \mkern -12mu \sqcap$

{\bf Proof of Theorem \ref{square2.theorem}.} We start as in the proof of Theorem \ref{square.theorem}
by considering the symmetrized process
$$P_n^{\epsilon} f^2 := \sum_{i=1}^n \epsilon_i f^2 (X_i)/n  $$
with $\epsilon_1 , \ldots , \epsilon_n$ a Rademacher sequence independent of $X_1 , \ldots , X_n $.
But when applying Dudley's theorem, see Theorem \ref{Dudley.theorem}, we use a different entropy bound.
For two functions $f $ and $\tilde f$ in the class ${\cal F}$
$$ \| f^2 - \tilde f^2  \|_n  \le
\| f+ \tilde f \|_{\infty} \| f - \tilde f \|_{n } \le2  K \| f - \tilde f \|_{n} .$$
It follows that
$${\cal H} ( u , {\cal F}^2 , \| \cdot \|_n ) \le {\cal H} ( u/ (2 K)  , {\cal F} , \| \cdot \|_{n} ) , \ u > 0 . $$
Hence
$$\int_{\delta}^{1} \sqrt {{\cal H} (u \hat R K , {\cal F}^2 , \| \cdot \|_n }du
 \le \int_{\delta}^{1} \sqrt {{\cal H} ( u\hat R /2   , {\cal F} , \| \cdot \|_{n} )} du. $$
So by Theorem \ref{Dudley.theorem}
$$\EE \biggl ( \sup_{f \in {\cal F} } P_n^{\epsilon} f^2 \biggr ) \le 
   K \EE {\cal J}_0 (\hat R, {\cal F} )/\sqrt n   .$$
   
Since $v \mapsto {\cal J}_0 ( \sqrt {v} , {\cal F}) $ is concave
$$ \EE {\cal J}_0 (\hat R, {\cal F} )   \le {\cal J}_0 ( \sqrt {\EE \hat R^2 } , {\cal F} ). $$
But then by Theorem \ref{symmetrization.theorem}
\begin {equation}\label{first-bound2.equation}
 \EE \biggl (\sup_{f \in {\cal F} } \biggl |  \| f \|_n^2 - \| f \|^2 \biggr | \biggr ) \le
2 K {\cal J}_0 (\sqrt { \EE \hat R^2} , {\cal F} ) / \sqrt {n}. 
\end{equation}
This leads to the by-product of the theorem:
$$ \EE \hat R^2 \le R^2 + 2K  {\cal J}_0 (\sqrt {\EE \hat R^2} , {\cal F} ) / \sqrt n
\le R^2 +  \EE \hat R^2/2  +  H (4 K /  \sqrt n ) . $$
or
$$\EE \hat R^2 \le 2 R^2 + H(4K / \sqrt n ) \le 4 R^2 . $$
Insert this back to find
$$ \EE \biggl (\sup_{f \in {\cal F} } \biggl |  \| f \|_n^2 - \| f \|^2 \biggr | \biggr ) \le 2 K {\cal J}_0  (2R
, {\cal F})  /\sqrt n  . $$ 
Finally apply Theorem \ref{bounded.theorem}. 

\hfill $\sqcup \mkern -12mu \sqcap$

{\bf Proof of Theorem \ref{Sobolev.theorem}.}

{\bf Case 1.} This follows from
$$\int_0^z \sqrt { H( u , {\cal F} , \| \cdot \|_{\infty})  } du \le \sqrt {p} \int_0^{z} 
 \sqrt{ H(u/p, {\cal F}_0 , \| \cdot \|_{\infty} ) } du \asymp \sqrt {p} p^{\alpha} z^{1 - \alpha} , $$
 and inserting this in Theorem \ref{square.theorem}.
 
 {\bf Case 2.} Here we use that by conditions (\ref{supnorm.equation}) and (\ref{incoherence.equation}), for $f = \sum_{k=1}^p f_{0,k} $,
 $$\| f \|_{\infty} \le c_2 \sum_{j=1}^p \| f_{0,k} \|^{1- \alpha} \le c_2 p^{1+ \alpha \over 2}\biggl (  \sum_{k=1}^p
 \| f_{0,k} \|^2 \biggr )^{1- \alpha  \over  2 } \le  c_2 p^{1+ \alpha \over 2} c_1^{1 - \alpha } 
 \| f \|^{1- \alpha} . $$
The result then follows applying the entropy bound of Case 1.

{\bf Case 3.}  For $f(x_1 , \ldots , x_p) = \sum_{k=1}^p \sum_{r=1}^N \beta_{r,k}  \psi_r (x_k) \in {\cal F}$ we have
$$\sum_{k=1}^p \sum_{r=1}^N  | \beta_{r,k} | \le p . $$ The result follows from Theorem \ref{l1.square.theorem}.

{\bf Case 4.} For $f(x_1, \ldots , x_p) = \sum_{k=1}^p \sum_{r=1}^N \beta_{r,k}  \psi_r (x_k) \in {\cal F}$, write
(with some abuse of notation)
$\beta_k:= ( \beta_{1,k} , \ldots , \beta_{N,k} )^T $ and $f_{\beta_k ,k} (x_1 , \ldots , x_p):=
\sum_{r=1}^N \beta_{r,k} ( x_k) $. Then by conditions (\ref{eigenvalue.equation}) and 
(\ref{incoherence.equation})
$$ \sum_{k=1}^p \| \beta_k \|^2 \le c_0  \sum_{k=1}^p N \| f_{\beta_k} \|^2  \le c_0 c_1N  \| f \|^2 . $$
The result then follows from Theorem \ref{l3.square.theorem}.

{\bf Case 5.} 
Let 
$$\tilde {\cal F} := \biggl \{  \tilde f (x_1 , \ldots , x_p) = \sum_{k=1}^p \sum_{r=1}^N \psi_r (x_k) \biggr \} . $$
For all $f \in {\cal F}(1)$ with $\| f \| \le 1 $ there is a $\tilde f \in \tilde {\cal F}$ such that
$$\| f - \tilde f \|_{\infty}  \le c_0 p N^{-{1 \over 2 \alpha } }  .$$
It follows that $\| \tilde f \| \le 1 + c_0 p N^{-{1 \over 2 \alpha } }  \le 2 $ for
$N \ge (c_0 p)^{2 \alpha}$. 
Define 
$${\bf Z}_0 ( \tilde {\cal F}(1)) := \sup_{\tilde f \in \tilde {\cal F}, \ \| \tilde f \| \le 1  } \biggl | \| \tilde f \|_n - \| \tilde f \| \biggr | . $$
Then
$$ {\bf Z}_0 (\tilde {\cal F}(1)) 
 = \sup_{ \tilde f \in \tilde {\cal F}, \ \| \tilde f \| \le 1 }
\biggl | {  \| \tilde f \|_n - \| \tilde f \|  \over \| \tilde f \| } \biggr |  \| \tilde f \| 
 \le \sup_{ \tilde f \in \tilde {\cal F}, \ \| \tilde f \| \le 1 }
\biggl | {  \| \tilde f \|_n - \| \tilde f \|  \over \| \tilde f \| } \biggr |   $$ $$
 = \sup_{ \tilde f \in \tilde {\cal F}, \ \| \tilde f \| = 1 }
\biggl | {  \| \tilde f \|_n - \| \tilde f \|   } \biggr |   
 = \sup_{ \tilde f \in \tilde {\cal F}, \ \| \tilde f \| = 1 }
\biggl | {  \| \tilde f \|_n^2 - \| \tilde f \|^2 \over \| \tilde f \|_n + \| \tilde f \|    } \biggr |   $$ $$
 \le \sup_{ \tilde f \in \tilde {\cal F}, \ \| \tilde f \| = 1 }
\biggl | {  \| \tilde f \|_n^2 - \| \tilde f \|^2 \over  \| \tilde f \|    } \biggr |   
 =  \sup_{ \tilde f \in \tilde {\cal F}, \ \| \tilde f \| = 1 }
\biggl | {  \| \tilde f \|_n^2 - \| \tilde f \|^2 \   } \biggr |   $$ $$ \le \sup_{ \tilde f \in \tilde {\cal F}, \ \| \tilde f \| \le 1 }
\biggl | {  \| \tilde f \|_n^2 - \| \tilde f \|^2 \   } \biggr |   
 := {\bf Z}( \tilde {\cal F}(1))= {\mathcal O} 
\biggl ( N \sqrt { p \log {(pN)} \over n} \biggr ) $$
where the last step follows from the same arguments as for Case 4. 
Define  now
$${\bf Z}_0 ( {\cal F}(1)) := \sup_{f \in  {\cal F}, \ \|  f \| \le 1  } \biggl | \| f \|_n - \| f \| \biggr | . $$
Clearly
$$ {\bf Z}_0 ( {\cal F}(1) )\le  2 {\bf Z}_0 ( \tilde {\cal F}(1) ) + 2 c_0 p N^{-{1 \over 2 \alpha} } . $$ 
So we find
$${\bf Z}_0 ( {\cal F}(1)) ={\mathcal O}_{\PP} (1)  \underbrace { \biggl ( N \sqrt { p \log {(pN)} \over n} \biggr )}_{I} +\underbrace{2 c_0 p N^{-{1 \over 2 \alpha} }}_{II}. $$

Since $p \le n$ choosing
$N\asymp \biggl ( { n p / \log (n)}  \biggr )^{\alpha \over  (1+ 2 \alpha) } $
gives
$$I \asymp II   \asymp p^{ {1 + 4 \alpha \over  2(1+2 \alpha)} }  \biggl ( { \log n / n} \biggr )^{1 \over
2 (1+ 2 \alpha) } ,$$
so that
${\bf Z}_0 ( {\cal F}(1)) = {\mathcal O}_{\PP}  ( p^{ {1 + 4 \alpha \over  2(1+2 \alpha)} }   ( { \log n / n} )^{1 \over
2 (1+ 2 \alpha) }  ) $.
But then also
$${\bf Z} ({\cal F}(1)) = \sup_{f \in {\cal F}, \ \| f \| \le 1 }
\biggl | \| f \|_n^2 - \| f \|^2 \biggr | $$ $$ \le 2  \sup_{f \in {\cal F}, \ \| f \| \le 1 }
\biggl | \| f \|_n - \| f \| \biggr |  +  \sup_{f \in {\cal F}, \ \| f \| \le 1 } \biggl | \| f \|_n - \| f \| \biggr |^2  $$
$$ = 2 {\bf Z}_0 ({\cal F} (1))+ {\bf Z}_0^2 ({\cal F} (1))= 
{\mathcal O}_{\PP} ( p^{ {1 + 4 \alpha \over  2(1+2 \alpha)} }  ( { \log n / n} )^{1 \over
2 (1+ 2 \alpha) }  ) . $$

\hfill $\sqcup \mkern -12mu \sqcap$

 {\bf Proof of Lemma \ref{Agmon.lemma}.}
 Let $N \in \Nat$ be arbitrary and let for 
$f_0 \in {\cal F}_0$ and certain $\beta_r$ and $\psi_r$, $r=1 , \ldots N$:
$$\| f_0 - \sum_{r=1}^N \beta_r \psi_r \|_{\infty} \le c_0 N^{-{1 \over 2\alpha}} . $$
Let $\beta := ( \beta_1 , \ldots  , \beta_N )^T$, $f_{\beta} := \sum_{r=1}^N \beta_r \psi_r $ and
$f_{\beta,k} (x_1 , \ldots , x_p) := f_{\beta} (x_k) $ and let
$f_{0,k} (x_1 , \ldots , x_p) := f_0 (x_k) $.
Then for $N \ge \|  f_{0,k} \|^{-2 \alpha  } $
$$ \| f_{\beta, k} \| \le \| f_{0,k} \| + c_0 N^{-{1 \over 2\alpha}} \le (1+ c_0) \| f_{0,k} \| . $$
Define $K_{\psi} := \sup_{r} \| \psi_r \|_{\infty} \vee 1$. We find for $N \ge 1$,
inserting (\ref{nonoverlap.equation}), 
$$\| f _0 \|_{\infty} \le\|  \sum_{r=1}^N \beta_r \psi_r \|_{\infty} + c_0 N^{-{1 \over 2 \alpha}} $$
$$ \le (1+c_3) \| \beta \|_{\infty} K_{\psi} +  c_0 N^{-{1 \over 2\alpha}} \le (1+c_3)\| \beta \|_2 K_{\psi}  + c_0 \| f_{0,k} \| $$
$$ \le c_0 (1+c_3) K_{\psi} \sqrt {N }\| f_{\beta, k} \| + c_0 \| f_{0,k} \| \le c_0 (2+ c_3 ) K_{\psi} \sqrt {N} \| f_{0,k} \|  $$
where in the second last inequality we used (\ref{eigenvalue.equation}). Take $N$ as the smallest integer
greater than or equal to $\|  f_{0,k} \|^{-2 \alpha  }$. Then $N \le \|  f_{0,k} \|^{-2\alpha  } +1 $
so that $ \sqrt {N}  \le \|  f_{0,k} \|^{-\alpha  } +1 $ and hence
$$ \| f _0 \|_{\infty} \le c_0 (2+ c_3 ) K_{\psi}(  \| f_{0,k} \|^{1- \alpha} + \| f_{0,k } \| ) \le 
2c_0(2+ c_3 )^2 K_{\psi}  \| f_{0,k} \|^{1- \alpha}, $$
since by (\ref{incoherence.equation}) $\| f_{0,k} \| \le 1 $.  Hence 
(\ref{supnorm.equation}) holds with $c_2 = 2c_0(2+ c_3 )^2 K_{\psi}$. 
\hfill $\sqcup \mkern -12mu \sqcap$

\subsection{Proofs for Section \ref{empirical.innerproduct.section}}

%{\bf Proof of Theorem \ref{product.theorem}.} 
%%We consider again first the symmetrized process
%%$$P_n^{\epsilon} fg  := \sum_{i=1}^n \epsilon_i f (X_i) g(X_i) , $$
%%with $\epsilon_1 , \ldots , \epsilon_n$ a Rademacher sequence independent of $X_1 , \ldots , X_n $.
%For functions $f , \tilde f $ and $g,\tilde g$ in the class ${\cal F}$ 
%$$ \| fg - \tilde f \tilde g \|_n \le \| f g - \tilde f g \|_n + \| \tilde f g - \tilde f \tilde g \|_n \le
%\hat R \| f - \tilde f \|_{n, \infty} + \hat R \| g - \tilde g \|_{n, \infty} $$
%with $\hat R := \sup_{f \in {\cal F}} \| f \|_n $.
%It follows that
%$${\cal H} ( u , {\cal F} \times {\cal F} , \| \cdot \|_n )   \le 2 {\cal H} ( u/ (2 \hat R)  , {\cal F} , \| \cdot \|_{n, {\infty} } ) , \ u > 0 . $$
%%Hence
%%$$\int_{\delta}^{1} \sqrt {{\cal H} (u \hat R K , {\cal F} \times {\cal F}  , \| \cdot \|_n }) du 
%% \le 2 \int_{\delta}^{1} \sqrt {{\cal H} ( uK /2   , {\cal F} , \| \cdot \|_{n, {\infty} } )} du .$$
%%So by Theorem \ref{Dudley.theorem}
%%$$\EE \biggl ( \sup_{f , g \in {\cal F} } P_n^{\epsilon} fg \biggr ) \le 
%% C_0 \inf_{\delta >0}  \EE \biggl [2 \hat R K \int_{\delta}^{1} \sqrt { {\cal H} (uK/2 , {\cal F} , \| \cdot \|_{n, {\infty} } )} du
%%/ \sqrt {n} 
%%+ \delta \hat R K   \biggr ]  $$
%%$$ \le  2 J_{\infty}  (K, {\cal F} ) \sqrt { \EE \hat R^2 }  /\sqrt n . $$
%We can further use the same arguments as in the proof of  Theorem \ref{square.theorem}.
%\hfill  $\sqcup \mkern -12mu \sqcap$

{\bf Proof of Theorem \ref{product2.theorem}.}
Let
$$ \hat R_1 :=  \sup_{f \in {\cal F} }\| f\|_n, \ \hat R_2 :=  \sup_{g \in {\cal G} }\| g\|_n. $$
For functions $f , \tilde f $ in the class ${\cal F}$ and $g,\tilde g$ in the class ${\cal G}$ we  have
$$ \| fg - \tilde f \tilde g \|_n \le \| f g - \tilde f g \|_n + \| \tilde f g - \tilde f \tilde g \|_n \le
\hat R_2 \| f - \tilde f \|_{n, \infty} + \hat R_1 \| g - \tilde g \|_{n, \infty} . $$
It follows that
$${\cal H} ( u , {\cal F} \times {\cal G} , \| \cdot \|_n )  $$ $$ \le {\cal H} ( u/ (2 \hat R_2)  , {\cal F} , \| \cdot \|_{n, {\infty} } ) +
{\cal H} ( u/ (2 \hat R_1)  , {\cal G} , \| \cdot \|_{n, {\infty} } ) , \ u > 0 . $$
We moreover have
$$\| fg \|_n \le(  \hat R_2 K_1) \wedge (\hat R_1 K_2) \le \hat R_1 K_2   .$$

Define the set
$${\cal R} := \{ \hat R_1 \le 2 R_1 , \  \hat R_2 \le  2 R_2  \} . $$
By Theorem \ref{square.theorem},  and since $J_{\infty} \le {\cal J}_{\infty}$, for values of $t$ and $n$ satisfying
 (\ref{R1.equation}) and (\ref{R2.equation})
  it holds that
  $$\PP ( {\cal R} ) \ge  1- 2\exp[-t] . $$

On ${\cal R}$
$$\int_{\delta \hat R_1 K_2}^{\hat R_1 K_2 } \sqrt {{\cal H} ( u , {\cal F} \times {\cal G} , \| \cdot \|_n ) } du \le $$
%$$2 R_1 \int_{\delta K_2}^{K_2} \sqrt {{\cal H} ( u/2 , {\cal G}  , \| \cdot \|_n ) } du 
% +  \int_{\delta R_1 K_2 /4}^{ 2   R_1  K_2  } \sqrt {{\cal H} ( u /(2 R_2) , {\cal F}  , \| \cdot \|_n ) } du $$
$$ = 2 R_1 \int_{\delta K_2}^{K_2 } \sqrt {{\cal H} ( u/2 , {\cal G}  , \| \cdot \|_n ) } du 
+ 2 R_2 \int_{\delta R_1 K_2 / (4R_2) }^{   K_2  R_1 / R_2  } \sqrt {{\cal H} ( u /2 , {\cal F}  , \| \cdot \|_n ) } du .$$

Consider now the symmetrized process
$$P_n^{\epsilon} fg  := \sum_{i=1}^n \epsilon_i f (X_i) g(X_i)/n  , $$
with $\epsilon_1 , \ldots , \epsilon_n$ a Rademacher sequence independent of $X_1 , \ldots , X_n $.
By Theorem \ref{Dudley.theorem} we have now found that conditionally
on ${\bf X} := \{ X_1 , \ldots, X_n \}$,
$$ \EE\biggl (  \sup_{f \in {\cal F}, \ g \in {\cal G} } | P_n^{\epsilon} fg   | \biggr \vert {\bf X}   \biggr ) 
\{ {\bf X} \in {\cal R} \} $$
%$$ 2 C_4 \EE \biggl [  \hat R_2  \int_{1/ ( 2 C_5 \sqrt n \hat R_2 )}^{K_1 } 
%\sqrt { {\cal H} ( u  , {\cal F} , \| \cdot \|_{n, {\infty} } ) } du \biggr ] / \sqrt n $$
%$$ + 2 C_4 \EE \biggl [  \hat R_1 \int_{1/ ( 2 C_5 \sqrt n \hat R_1)}^{K_2}  
%\sqrt { {\cal H} ( u  , {\cal G} , \| \cdot \|_{n, {\infty} } )} du \biggr ] / \sqrt n  . $$
$$ \le 2 R_1 {\cal J}_{\infty} ( K_2 , {\cal G} )  / \sqrt n + 2 R_2 {\cal J}_{\infty} ( R_1 K_2  / R_2, {\cal F} ) / \sqrt n := {\bf E}  $$
From Theorem \ref{bounded.theorem}, we get that conditionally on ${\bf X} $
$$ \PP \biggl (  \sup_{f \in {\cal F}, \ g \in {\cal G} } | P_n^{\epsilon} fg   |/ C_1  \ge
{\bf E}  + 2 R_1 K_2 \sqrt {t \over n} + K_1 K_2 { t \over n } \biggr \vert {\bf X}  
 \biggr ) \{ {\bf X} \in {\cal R} \} \le \exp [-t] . $$
 But then, since $\PP ({\cal R}^c ) \le 2 \exp [-t] $, 
 $$  \PP \biggl (  \sup_{f \in {\cal F}, \ g \in {\cal G} } | P_n^{\epsilon} fg   |/ C_1  \ge
 {\bf E} + 2 R_1 K_2 \sqrt {t \over n} + K_1 K_2 { t \over n } 
 \biggr )  \le 3 \exp [-t] . $$
 Now apply Theorem \ref{symmetrization2.theorem}. 
  \hfill $\sqcup \mkern -12mu \sqcap$ 
  
  {\bf Proof of Theorem \ref{subGaussianproduct.theorem}.}  Let $\{ \varepsilon_i \}_{i=1}^n $ be a Rademacher sequence independent of
$\{ X_i, Y_i \}_{i=1}^n $. 
Conditionally on $( {\bf X} , {\bf Y}) := ( \{ X_1 , \ldots , X_n \}, \{ Y_1 , \ldots , Y_n \}) $, by  Theorem \ref{Dudley.theorem}, for $P_n^{\varepsilon}  {\bf Y} f := \sum_{i=1}^n \varepsilon_i Y_i f(X_i) / n $, $f \in {\cal F}$:
$$\EE \biggl ( \sup_{f \in {\cal F} } | P_n^{\varepsilon} {\bf Y} f| \biggr  \vert {\bf  X} , {\bf Y} \biggr ) 
% $$ $$ \le 
% C_0   \biggl [ \hat R  \int_{0 }^{C_3  \hat R K_0 } \sqrt { {\cal H} (u / (C_3 K_0) , {\cal F} , \| \cdot \|_n  )} du \biggr ]  / \sqrt {n}  $$ $$
  \le C_0  {\cal J}_{\infty}  ( \| {\bf Y } \|_n K , \{ Y f : \ f \in {\cal F}\}  ) . $$
 So on the set
 $${\cal Y} := \{ \| {\bf Y} \|_n^2 \le 4 K_0^2 \} $$
 we get
$$\EE \biggl ( \sup_{f \in {\cal F} } | P_n^{\varepsilon} {\bf Y } f | \biggr  \vert {\bf  X}, {\bf Y}  \biggr )\{  
{\bf Y} \in {\cal Y}  \} 
% $$ $$ \le 
% 2 C_0 \inf_{\delta > 0 } \biggl [  C_3 K  \int_{ \delta  R  /4 }^{  R } \sqrt { {\cal H}_{n, {\infty} } (u/2 , {\cal F}   )} du \biggr ]  / \sqrt {n} +  C_3  R K_0  \delta  \biggr  ] $$ $$ 
 \le  2 {\cal J }_{\infty}  (K_0 K , {\cal F} )   / \sqrt n .$$
Now apply Theorem \ref{subGaussian.theorem} to obtain that for all $t >0$
$$\PP \biggl (\sup_{f \in {\cal F}}  | P_n^{\varepsilon} {\bf Y}  f  | /C_1  \ge
{   {\cal J}_{\infty}  (K K_0 , {\cal F} )  \over  \sqrt n  }  + 2 K K_0R  \sqrt {t \over n}  \biggr \vert {\bf X} , {\bf Y}  \biggr )  \{ {\bf Y} \in {\cal Y} \} 
  \le \exp[-t] . $$
  We now integrate out and use that   $\PP ({\cal Y}^c ) \le \exp [-t] $.
  Then we de-symmetrize using
  Theorem \ref{symmetrization2.theorem}. 
  
\hfill $\sqcup \mkern -12mu \sqcap$

\subsection{Proofs for Section \ref{linear.section}} 

{\bf Proof of Theorem \ref{l1.square.theorem}.}

If $\| \beta \|_1 \le M $ we know that 
$$\| f_{\beta} \|_{\infty} \le \|\beta \|_1 
\max_{1 \le i \le n } \max_{1 \le j \le p } |X_{i,j} | \le MK_X . $$
Fixing $\delta $ at $\delta = 1/\sqrt n $ in (\ref{J.definition}), we find by Lemma \ref{Rudelson.lemma}
$$J({\cal F}) \le C_0 MK_X
\int_{1/  \sqrt {n}}^{1} \sqrt { {\cal H} (u/2 , \{f_{\beta} : \  \| \beta \|_1 \le M \} , \| \cdot \|_{n, {\infty} } )} du +
C_0 MK $$ 
$$=2C_0  \sqrt {\log (2p)\log (2n)} MK_X \log n + C_0 MK_X . $$
The result now follows from Theorem \ref{square.theorem}.
\hfill $\sqcup \mkern -12mu \sqcap$

{\bf Proof of Theorem \ref{l0.square.theorem}.} 
Since $\| f_{\beta_S} \| \le 1 $ implies
$$\| \beta_S \|_1 \le \sqrt {s} \| \beta \|_2 \le \sqrt {s} / \Lambda_{\rm min}  $$
this follows directly from
Theorem \ref{l1.square.theorem}. 
\hfill $\sqcup \mkern -12mu \sqcap$

{\bf Proof of Lemma \ref {Gaus.lemma}.}
We let $\psi_j (X_i) := X_{i, j} $, $i=1 , \ldots , n$, $j=1 , \ldots , p$. 
It holds that
$$\sup_{\| f_{\beta} \| \le 1 }|  ( P_n - P){\bf Y} f_{\beta} |^2 \le \sup_{\| f_{\beta} \| \le 1 }
\sum_{j=1}^p | (P_n - P) {\bf Y} \psi_j |^2 \| \beta \|_2^2 \le
\sum_{j=1}^p | (P_n - P) {\bf Y } \psi_j |^2 / \Lambda_{\rm min}^2 . $$
But by the triangle inequality
$$\biggl \| \sum_{j=1}^p | (P_n - P) {\bf Y} \psi_j|^2 \biggr  \|_{\Psi_1}  \le \sum_{j=1}^p\biggl \| | (P_n - P) {\bf Y} \psi_j|^2\biggr \|_{\Psi_1 } . $$
Moreover for all $j$,
$$ \biggl  \| | (P_n - P) {\bf Y}  \psi_j|^2 \biggr  \|_{\Psi_1 } = \biggl \| (P_n - P)  {\bf Y} \psi_j \biggr \|_{\Psi_2}^2 \le c_2^2  K_0^2 K_X^2 /n . $$
Hence
$$\biggl \| \sqrt {\sum_{j=1}^p | (P_n - P) {\bf Y}  \psi_j|^2 } \biggr \|_{\Psi_2}= \sqrt {\biggl \| \sum_{j=1}^p | (P_n - P) {\bf Y}  \psi_j|^2  
\biggr \|_{\Psi_1}}
\le  c_2 K_0 K_X \sqrt {p \over n }.  $$
By Chebyshev's inequality, for all $t > 0 $, 
$$\PP \biggl ( \sqrt {\sum_{j=1}^p | (P_n - P) {\bf Y}  \psi_j |^2   } \ge { c_2 K_0 K_X \over \Lambda_{\rm min} }
\sqrt { p t \over n }  \biggr ) \le 2 \exp[-t] . $$
%But
%$$\EE \sum_{j=1}^p | (P_n - P) Y \psi_j |^2 = \sum_{j=1}^p P(Y \psi_j)^2 / n \le
%p K_0^2 K_X^2 . $$
%Hence
%$$\EE \sup_{\| f_{\beta} \| \le 1 }|  ( P_n - P)Y f_{\beta} | \le \sqrt {p} K_0 K_X / \Lambda_{\rm min}. $$
%The result follows from Theorem \ref{sub-exponential.theorem}. 
\hfill $\sqcup \mkern -12mu \sqcap$

\subsection{Proofs for Section \ref{least.squares.section}} 

{\bf Proof of Lemma \ref{basic.lemma}.}
The inequality
$$\| {\bf Y} - \hat f \|_n \le \| {\bf Y}  - f^* \|_n , $$
can be rewritten to the Basic Inequality
$$\| \hat f - f^* \|_n^2 \le 2 P_n ({\bf Y}- f^*)  ( \hat f - f^*) . $$
On ${\cal T}$ we therefore have
$$ \| \hat f - f^* \|^2 \le \| \hat f - f^* \|^2 - \| \hat f - f^* \|_n^2 + 2 P_n ({\bf Y}-f^*) (\hat f - f^*) $$
$$ \le \delta_n \|  \hat f - f^* \|^2 + \delta_n \| \hat f - f^* \| , $$
where we used that $P ({\bf Y}-f^*) (\hat f - f^*)=0 $. 
Hence
$$\| \hat f - f^* \| \le \delta_n / (1- \delta_n ) \le 2 \delta_n  . $$

\hfill $\sqcup \mkern -12mu \sqcap$

{\bf Proof of Theorem \ref{theoretical.norm.theorem}.}  
This follows from Lemma \ref {basic.lemma} combined with Summary \ref{linear.summary}. 
We use here that $\| f^* \| \le \|{\bf Y} \| \le K_0 $.
The first result then follows immediately from Lemma \ref {basic.lemma} and Summary \ref{linear.summary}.
For the second result, write
$$\| {\bf Y} - \hat f \|_n^2 - \| {\bf Y}  - f^* \|_n^2 = \| \hat f - f^* \|_n^2 - 2 P_n ({\bf Y} - f^*) (\hat f - f^*) .$$
But
$$\| \hat f - f^* \|_n^2 = \| \hat f - f^* \|^2 + \biggl ( { \| \hat f - f^* \|_n^2 \over \| \hat f - f^* \|^2}  -1 \biggr ) 
\| \hat f - f^* \|^2 $$
$$ = {\mathcal O}_{\PP} (\delta_n^2 + \delta_n^3 ) $$
and
$$  P_n  ({\bf Y}- f^*) (\hat f - f^*) = (P_n-P)   ({\bf Y} - f^*) (\hat f - f^*) $$ $$= \| \hat f - f^* \| \biggl ( (P_n-P)   ({\bf Y}- f^*) (\hat f - f^*)/ \| \hat f - f^* \| \biggr ) = O_{\PP} (\delta_n^2 ) . $$
Hence
$$\| {\bf Y} - \hat f \|_n^2 - \| {\bf Y} - f^* \|^2 = \| {\bf Y} - \hat f \|_n^2 - \| {\bf Y} - f^*\|_n^2 + \| {\bf Y} - f^* \|_n^2 +
\| {\bf Y} - f^* \|_n^2  $$ $$= O_{\PP} (\delta_n^2 ) +  \|{\bf  Y} - f^* \|_n^2 -\| {\bf Y} - f^* \|^2 . $$
We find
$$ \| {\bf Y} - f^* \|_n^2 -\| {\bf Y} - f^* \|^2 = \| {\bf Y} \|_n^2 - \| {\bf Y} \|^2 - 2 (P_n - P) ( {\bf Y} f^* ) + \| f^* \|_n^2 -
\| f^* \|^2$$ $$ = {\mathcal O}_{\PP} ( 1/ \sqrt n ) + {\mathcal O} _{\PP} (\delta_n) .$$
   \hfill $\sqcup \mkern -12mu \sqcap$

{\bf Proof of Theorem 
 \ref{theoretical.norm2.theorem}.}   
 By Summary \ref{linear.summary} all probability statements
are uniformly in $S$, so that the set ${\cal T}$ given in Lemma \ref{basic.lemma} has with
\break $\delta_n = O(\sqrt { p  \log p / n })$ the required large probability. 
    
    \hfill $\sqcup \mkern -12mu \sqcap$
    
    \subsection{Proof for Section \ref{DAG.section}.}
    
    {\bf Proof of Theorem \ref{DAG0.theorem}.} 
%Let
%$${\cal F} := \biggl \{ f = \sum_{j=1}^p f_{0,j} , \ f_{0,j} (x_1 , \ldots , x_p) = f_0 (x_j) , \ f_0 \in {\cal F}_0 \biggr \} . $$
%Then from Theorem \ref{square.theorem}, we have for $p^3 / n\rightarrow 0 $
From Theorem \ref{Sobolev.theorem}
$$\sup_{f \in {\cal F} , \ \| f \| \le 1 } 
\biggl | \| f \|_n^2 - \| f \|^2\biggr | = o_{\PP } (1) . $$
We know moreover from Lemma \ref{subGaussianproduct.theorem} that also for 
 $$\max_j \sup_{f \in {\cal F} , \ \| f \| \le 1 } ( P_n -P) {\bf X}_j f = o_{\PP} (1) .$$
 Also
 $$\max_{j} \biggl | \| {\bf X}_j \|_n^2   - \| X_j \|^2 \biggr | = o_{\PP} (1) . $$
 Hence,
 \begin{equation} \label{uniform.equation}
 \max_j \sup_{f \in {\cal F}, \ \| f \| \le 1  } \biggl | \| {\bf X}_j - f \|_n^2 - \|{ \bf X}_j - f \|^2  \biggr | = o_{\PP} (1) . 
 \end{equation}

%If we assume in addition that 
%$$\| f_{0,j} \|_{\infty} \le c_1 \| f_{0,j} \|^{1- \alpha}, f_{0,j} (x_1 , \ldots , x_p ) = f_0 (x_j) , \ f_0 \in {\cal F}_0  $$
%and
% $$ \sum_{j=1}^p \| f_{0,j} \|^2 \le c_2 \| \sum_{j=1}^p f_{0,j} \|^2 , f_{0,j} (x_1 , \ldots , x_p ) = f_0 (x_j) , \ f_0 \in {\cal F}_0  $$
% where $\rho= {\mathcal O}(1)$, then we get the uniform convergence 
% (\ref{uniform.equation}) for $p^{3-(1- \alpha)^2 } / n \rightarrow 0 $.
% 
 We now note that we only need uniform convergence over $f \in {\cal F}$ with $\| f \| \le 1$. To see
 this, let for any $\pi$ and $j$
 $$\tilde f_j (\pi) := s \hat f_j (\pi) + (1-s) f_j^* (\pi)  $$
 where $s:= c/ (c+ \| \hat f_j (\pi) - f_j^* (\pi) \|) $ and $c$ a constant (depending on $j$ and $\pi$)
 to be chosen (see below).
 Then $\| \tilde f_j (\pi) - f_j^* ( \pi ) \| \le 1 $. Moreover
 $$\| {\bf X}_j - \tilde f_j (\pi) \|_n^2  \le s \| {\bf X}_j - \hat f_j (\pi) \|_n^2 + (1-s) \|
 {\bf X}_j - f_j^* (\pi) \|_n^2 \le \| {\bf X}_j - f_j^* ( \pi) \|_n^2 . $$
 So
 $$\| {\bf X}_j - \tilde f_j (\pi) \| = \| {\bf X}_j - \tilde f_j (\pi) \| -  \| {\bf X}_j - \tilde f_j (\pi) \|_n +
 \| {\bf X}_j - \tilde f_j (\pi) \|_n  $$
 $$ = \| {\bf X}_j - \tilde f_j (\pi) \|_n + o_{\PP} (1) \le \| {\bf X}_j -  f_j^* (\pi) \|_n + o_{\PP} (1) $$ 
 $$= \| {\bf X}_j -  f_j^* (\pi) \| + o_{\PP} (1) \le \| {\bf X}_j - \tilde f_j (\pi) \| +o_{\PP} (1) $$
 where in the last step we used the convexity of ${\cal F}_0$. 
 Thus $\| {\bf X}_j - \tilde f_j (\pi) \| = o_{\PP} (1) $. This impies that $\| \tilde f_j (\pi) - f_j^* (\pi) \|
 \le \| {\bf X}_j - f_j^* (\pi)  \| + o_{\PP} (1) $. Choosing $c$ appropriately, for example
 $c= 4 \| {\bf X}_j - f_j^* (\pi) \| $, we now find that $\| \hat f_j (\pi) - f_j^* (\pi) \|= {\mathcal O}_{\PP} (1)$.
 By applying the same arguments as above with $\tilde f_j (\pi)$ replaced by $\hat f_j (\pi)$
shows that $\| {\bf X}_j - \hat f_j (\pi) \| = o_{\PP} (1)$. This result is uniformly in $\pi \in \Pi$ by the same arguments
as used for  Theorem \ref{theoretical.norm2.theorem}. Application of the
union bound and deviation bounds for each $j$, we see that the result is also uniformly in $j $.

 \hfill $\sqcup \mkern  -12mu \sqcap$

\bibliographystyle{plainnat}
\bibliography{reference1}

\begin{thebibliography}{28}
\providecommand{\natexlab}[1]{#1}
\providecommand{\url}[1]{\texttt{#1}}
\expandafter\ifx\csname urlstyle\endcsname\relax
  \providecommand{\doi}[1]{doi: #1}\else
  \providecommand{\doi}{doi: \begingroup \urlstyle{rm}\Url}\fi

\bibitem[Adamczak et~al.(2011)Adamczak, Litvak, Pajor, and
  Tomczak-Jaegermann]{adamczak2011sharp}
R.~Adamczak, A.~E Litvak, A.~Pajor, and N.~Tomczak-Jaegermann.
\newblock Sharp bounds on the rate of convergence of the empirical covariance
  matrix.
\newblock \emph{Comptes Rendus Mathematique}, 349\penalty0 (3):\penalty0
  195--200, 2011.

\bibitem[Agmon and Jones(1965)]{Agmon:65}
S.~Agmon and F.~Jones.
\newblock \emph{{Lectures on elliptic boundary value problems Elliptic boundary
  value problems Van Nostrand mathematical studies}}.
\newblock Van Nostrand, 1965.

\bibitem[Ahlswede and Winter(2002)]{ahlswede2002strong}
R.~Ahlswede and A.~Winter.
\newblock Strong converse for identification via quantum channels.
\newblock \emph{Information Theory, IEEE Transactions on}, 48\penalty0
  (3):\penalty0 569--579, 2002.

\bibitem[Bartlett et~al.(2012)Bartlett, Mendelson, and Neeman]{bartlett2012}
P.L. Bartlett, S.~Mendelson, and J.~Neeman.
\newblock $\ell_1$-regularized linear regression: persistence and oracle
  inequalities.
\newblock \emph{Probability {T}heory and {R}elated {F}ields}, 154\penalty0
  (1-2):\penalty0 193--224, 2012.

\bibitem[Bickel et~al.(2009)Bickel, Ritov, and Tsybakov]{bickel2009sal}
P.~Bickel, Y.~Ritov, and A.~Tsybakov.
\newblock Simultaneous analysis of {L}asso and {D}antzig selector.
\newblock \emph{Annals of Statistics}, 37:\penalty0 1705--1732, 2009.

\bibitem[B\"uhlmann et~al.(2013)B\"uhlmann, Peters, and Ernest]{PetersB}
P.~B\"uhlmann, J.~Peters, and J.~Ernest.
\newblock {CAM: Causal Additive Models, high-dimensional order search and
  penalized regression}, 2013.
\newblock ArXiv 1310.1533.

\bibitem[Chen et~al.(1998)Chen, Donoho, and Saunders]{chen1998atomic}
S.S. Chen, D.L. Donoho, and M.A. Saunders.
\newblock Atomic decomposition by basis pursuit.
\newblock \emph{SIAM journal on scientific computing}, 20\penalty0
  (1):\penalty0 33--61, 1998.

\bibitem[Dudley(1967)]{dudley1967sizes}
R.M. Dudley.
\newblock The sizes of compact subsets of hilbert space and continuity of
  gaussian processes.
\newblock \emph{Journal of Functional Analysis}, 1\penalty0 (3):\penalty0
  290--330, 1967.

\bibitem[Gin{\'e} and Koltchinskii(2006)]{gine2006concentration}
E.~Gin{\'e} and V.~Koltchinskii.
\newblock Concentration inequalities and asymptotic results for ratio type
  empirical processes.
\newblock \emph{The Annals of Probability}, 34\penalty0 (3):\penalty0
  1143--1216, 2006.

\bibitem[Gu{\'e}don et~al.(2007)Gu{\'e}don, Mendelson, Pajor, and
  Tomczak-Jaegermann]{guedon2007subspaces}
O.~Gu{\'e}don, S.~Mendelson, A.~Pajor, and N.~Tomczak-Jaegermann.
\newblock Subspaces and orthogonal decompositions generated by bounded
  orthogonal systems.
\newblock \emph{Positivity}, 11\penalty0 (2):\penalty0 269--283, 2007.

\bibitem[Koltchinskii(2013)]{koltchinskii2011remark}
V.~Koltchinskii.
\newblock A remark on low rank matrix recovery and noncommutative {B}ernstein
  type inequalities.
\newblock In \emph{IMS Collections From Probability to Statistics and Back:
  High-Dimensional Models and Processes}, volume~9, pages 213--226. Institute
  of Mathematical Statistics, Beachwood, Ohio, 2013.
\newblock Banerjee, M., Bunea, F., Huang, J., Koltchinskii, V., and Maathuis,
  M. H., eds.

\bibitem[Ledoux and Talagrand(1991)]{Ledoux:91}
M.~Ledoux and M.~Talagrand.
\newblock \emph{{Probability in Banach Spaces: Isoperimetry and Processes}}.
\newblock Springer Verlag, New York, 1991.

\bibitem[Loh and Wainwright(2012)]{loh2012}
P.-L. Loh and M.J. Wainwright.
\newblock High-dimensional regression with noisy and missing data: {P}rovable
  guarantees with non-convexity.
\newblock \emph{Annals of Statistics}, 40:\penalty0 1637--1664, 2012.

\bibitem[Massart(2000)]{Massart:00a}
P.~Massart.
\newblock {About the constants in Talagrand's concentration inequalities for
  empirical processes}.
\newblock \emph{Annals of Probability}, 28:\penalty0 863--884, 2000.

\bibitem[Meier et~al.(2009)Meier, Van~de Geer, and B{\"u}hlmann]{meier2009high}
L.~Meier, S.~Van~de Geer, and P.~B{\"u}hlmann.
\newblock High-dimensional additive modeling.
\newblock \emph{The Annals of Statistics}, 37\penalty0 (6B):\penalty0
  3779--3821, 2009.

\bibitem[M\"uller and van~de Geer(2013)]{MullervdGeer13}
P.~M\"uller and S.A. van~de Geer.
\newblock The partial linear model in high dimensions, 2013.
\newblock Submitted, arXiv:1307.1067.

\bibitem[Nickl and {van de Geer}(2013)]{Nickl:2013}
R.~Nickl and S.A. {van de Geer}.
\newblock Confidence sets in sparse regression, 2013.
\newblock arXiv:1209.1508v2, to appear in The Annals of Statistics.

\bibitem[Pollard(1984)]{pollard1984convergence}
D.~Pollard.
\newblock \emph{{Convergence of Stochastic Processes}}.
\newblock Springer, 1984.

\bibitem[Raskutti et~al.(2010)Raskutti, Wainwright, and
  Yu]{raskutti2010restricted}
G.~Raskutti, M.J. Wainwright, and B.~Yu.
\newblock Restricted eigenvalue properties for correlated {G}aussian designs.
\newblock \emph{Journal of Machine Learning Research}, 11:\penalty0 2241--2259,
  2010.

\bibitem[Rudelson and Vershynin(2008)]{rudelson2008sparse}
M.~Rudelson and R.~Vershynin.
\newblock On sparse reconstruction from {F}ourier and {G}aussian measurements.
\newblock \emph{Communications on Pure and Applied Mathematics}, 61\penalty0
  (8):\penalty0 1025--1045, 2008.

\bibitem[Rudelson and Zhou(2013)]{rudelson2011reconstruction}
M.~Rudelson and S.~Zhou.
\newblock Reconstruction from anisotropic random measurements.
\newblock \emph{IEEE Transactions on Information Theory}, 59:\penalty0
  3434--3447, 2013.

\bibitem[Talagrand(1995)]{Talagrand:95}
M.~Talagrand.
\newblock {Concentration of measure and isoperimetric inequalities in product
  spaces}.
\newblock \emph{Publications Math\'ematiques de l'IHES}, 81:\penalty0 73--205,
  1995.

\bibitem[van~de Geer(2007)]{vandeG07}
S.~van~de Geer.
\newblock The deterministic {L}asso.
\newblock \emph{The JSM Proceedings}, 2007.

\bibitem[van~de Geer and B\"uhlmann(2009)]{vdG:2009}
S.~van~de Geer and P.~B\"uhlmann.
\newblock {On the conditions used to prove oracle results for the Lasso}.
\newblock \emph{Electronic Journal of Statistics}, pages 1360--1392, 2009.

\bibitem[{van de Geer} and B{\"u}hlmann(2013)]{vdGP:2013}
S.A. {van de Geer} and P.~B{\"u}hlmann.
\newblock $\ell_0$-penalized maximum likelihood for sparse directed acyclic
  graphs.
\newblock \emph{The Annals of Statistics}, 41:\penalty0 536--567, 2013.

\bibitem[van~de Geer and Mammen(2013)]{vdGMammen13}
S.A. van~de Geer and E.~Mammen.
\newblock Penalized least squares for an additive model, 2013.
\newblock in progress.

\bibitem[{van de Geer} et~al.(2013){van de Geer}, B{\"u}hlmann, and
  Ritov]{vdGBR:2013}
S.A. {van de Geer}, P.~B{\"u}hlmann, and Y.~Ritov.
\newblock On asymptotically optimal confidence regions and tests for
  high-dimensional models, 2013.
\newblock Submitted, arXiv:1303.0518.

\bibitem[{van der Vaart} and Wellner(1996)]{vanderVaart:96}
A.~W. {van der Vaart} and J.~A. Wellner.
\newblock \emph{{Weak Convergence and Empirical Processes}}.
\newblock Springer Series in Statistics. Springer-Verlag, New York, 1996.
\newblock ISBN 0-387-94640-3.

\end{thebibliography}

\end{document}